\documentclass[12pt, reqno]{amsart}
\usepackage{amssymb, latexsym, mathrsfs, pstcol, fullpage}
\newtheorem{thm}{Theorem}[section]
\newtheorem{corollary}[thm]{Corollary}
\newtheorem{lemma}[thm]{Lemma}
\newtheorem{theorem}[thm]{Theorem}

\newtheorem{prop}[thm]{Proposition}
\numberwithin{equation}{section}

\theoremstyle{definition}
\newtheorem{rem}[thm]{Remark}
\newtheorem{ex}[thm]{Example}

\newcommand{\al}{\alpha}
\renewcommand{\b}{\beta}

\newcommand{\de}{\delta}
\newcommand{\e}{\varepsilon}
\newcommand{\la}{\lambda}
\renewcommand{\phi}{\varphi}

\renewcommand{\d}{\partial}
\newcommand{\R}{{\mathbb R}}

\newcommand{\br}[1]{\left\langle #1 \right\rangle}
\newcommand{\Case}[1]{\noindent \underline{Case #1:}}
\newcommand{\Step}[1]{\noindent \underline{Step #1:}}
\newcommand{\Subcase}[1]{\noindent \underline{Subcase #1:}}

\renewcommand{\qed}{\rule{3mm}{3mm}}
\renewenvironment{proof}
    {\vspace{1mm}\noindent\textbf{Proof.}}
    {\hspace*{\fill} $\qed$\vspace{1mm}}
\newenvironment{proof_of}[1]
    {\vspace{1mm}\noindent {\bf Proof of #1.}}
    {\hspace*{\fill} $\qed$\vspace{1mm}}
\renewcommand{\qed}{\rule{3mm}{3mm}}

\begin{document}
\title[Nonlinear Waves with a Sign-Changing Potential]{Existence and Blow up of Small-Amplitude Nonlinear Waves with
A Sign-Changing Potential}
\author{Paschalis Karageorgis}
\address{Brown University, Box 1917, Providence, RI 02912}
\email{petekara@math.brown.edu}



\begin{abstract}
We study the nonlinear wave equation with a sign-changing potential in any space dimension.  If the potential is small
and rapidly decaying, then the existence of small-amplitude solutions is driven by the nonlinear term.  If the potential
induces growth in the linearized problem, however, solutions that start out small may blow-up in finite time.
\end{abstract}
\maketitle

\section{Introduction}
Consider the nonlinear wave equation with potential
\begin{equation}\label{we}
\left\{
\begin{array}{rll}
\d_t^2 u - \Delta u + V(x)\cdot u = F(u) \quad\quad &\text{in\: $\R^n \times (0,\infty)$} \\
u(x,0) = \phi(x) ;\quad \d_t u(x,0) = \psi(x) &\text{in\: $\R^n$,}
\end{array}\right.
\end{equation}
where $V(x)$ is some known function and $F(u)$ behaves like $|u|^p$ for some $p>1$.  When it comes to the special case
$V(x)\equiv 0$, this equation has been extensively studied since Fritz John's seminal work \cite{Jo}.  For that case, in
particular, the existence of small-amplitude solutions is known to depend on both the exact value of $p$ and the decay
rate of the initial data.  In this paper, we address the more general case \eqref{we} when the potential $V(x)$ is of
arbitrary sign. Our aim is to show that the existence of small-amplitude solutions may also be affected by two additional
parameters, namely, the amplitude and the decay rate of $V(x)$.

First, consider solutions to \eqref{we} when $V(x)\equiv 0$ and the small initial data have compact support. John's
classical result \cite{Jo} in $n=3$ space dimensions ensures their global existence if $p> 1+\sqrt 2$ and their blow-up
if $1< p < 1+\sqrt 2$.  More generally, a similar dichotomy holds in $n\geq 2$ space dimensions, where the borderline
case is given by the positive root $p_n$ of the quadratic
\begin{equation}\label{pn}
(n-1)p_n^2 = (n+1)p_n + 2;
\end{equation}
see \cite{GLS, Gl2, JZ, Jo, LS, Sid, Tat, YZ1}.  As for the borderline case $p= p_n$ with $n\geq 2$, the blow-up of
solutions persists \cite{Sc, YZ2}. Finally, when $n=1$, blow-up occurs for any $p>1$; see \cite{Kat}.

Next, consider solutions to \eqref{we} when $V(x)\equiv 0$ and the small initial data decay slowly.  In $n=2,3$ space
dimensions, their global existence is ensured as long as $p> p_n$ and the initial data satisfy
\begin{equation}\label{oec}
\sum_{|\al|\leq 3} \:|\d_x^\al \phi(x)| + \sum_{|\al|\leq 2} \:|\d_x^\al \psi(x)| \leq \e (1+|x|)^{-k-1}
\end{equation}
for some $k\geq 2/(p-1)$ and some small $\e>0$.  On the other hand, blow-up may occur for any $p>1$ when the initial data
are such that
\begin{equation}\label{obc}
\phi(x)=0, \quad \psi(x) \geq \e (1+|x|)^{-k-1} \quad\quad\text{in\, $\R^n$}
\end{equation}
for some $0\leq k< 2/(p-1)$ and $\e>0$; see \cite{AT, As, Ku, Ts1, Ts2, Ts3}.  In $n\geq 4$ space dimensions, the same
blow-up result holds, provided that $\phi,\psi$ are radially symmetric \cite{Ta1, Ta2}.  However, the existence result is
slightly modified as follows.  Instead of \eqref{oec}, one assumes that
\begin{equation}\label{oec2}
\sum_{|\al|\leq 2} \br{x}^{|\al|} \,|\d_x^\al \phi(x)| + \sum_{|\al|\leq 1} \br{x}^{|\al|+1} \,|\d_x^\al \psi (x)| \leq
\e \br{x}^{-k},
\end{equation}
where $\phi,\psi$ are radially symmetric and $\br{x}= 1+|x|$ for each $x\in \R^n$.  When $k\geq 2/(p-1)$ and $\e>0$ is
small, one then has global solutions in $n\geq 4$ space dimensions as well \cite{Kb1, KKe}, but the additional assumption
$2/(p-1) \neq k > n/2$ is imposed for even values of $n$.

In the remaining of this paper, we shall mostly focus on the radially symmetric version of the nonlinear wave equation
with potential \eqref{we}.  Thus, the equation of interest is
\begin{equation}\label{nl}
\left\{
\begin{array}{rll}
\d_t^2 u -\d_r^2 u -\dfrac{n-1}{r} \cdot \d_r u = F(u) - V(r)\cdot u \quad\quad &\text{in\, $\Omega_T= \R_+\times
(0,T)$}\\
u(r,0) = \phi(r); \quad \d_t u(r,0)= \psi(r) \quad\quad &\text{in\, $\R_+$.}
\end{array}\right.
\end{equation}
Before we state our main results, however, let us first introduce some hypotheses.  When it comes to the nonlinear term
$F(u)$, we shall impose the conditions
\begin{equation}\label{F}
F\in \mathcal{C}^1(\R); \quad F(0)= F'(0) = 0; \quad |F'(u)-F'(v)| \leq Ap |u-v|^{p-1}
\end{equation}
for some $A>0$ and some $p$ larger than the critical power $p_n$ \eqref{pn}.  When it comes to the potential term $V(r)$,
we require that
\begin{equation}\label{V}
\sum_{i=0}^1 \br{r}^i \,|V^{(i)}(r)| \leq V_0 \br{r}^{-\kappa}
\end{equation}
for some $V_0>0$ and $\kappa>2$.  As for the initial data, our exact assumption depends on the parity of $n$. In
particular, setting
\begin{equation}\label{m}
m = \left\{
\begin{array}{ccl}
(n-3)/2 &&\text{if\, $n$ is odd} \\
(n-2)/2 &&\text{if\, $n$ is even,}
\end{array}\right.
\end{equation}
we shall consider initial data $\phi,\psi$ such that
\begin{equation}\label{da3}
\sum_{i=0}^2 r^i \,|\phi^{(i)}(r)| + \sum_{i=0}^1 r^{i+1} \,|\psi^{(i)}(r)| \leq \e r^{1-m} \br{r}^{m-1-k}
\end{equation}
for some $\e>0$ and $k\geq 0$.  We remark that $m\geq 1$ when $n\geq 4$ and that \eqref{oec2} implies \eqref{da3} for
each $m\geq 1$.

The existence result of this paper can now be stated as follows.

\begin{theorem}\label{et}
Let $n\geq 4$ and define $m$ by \eqref{m}.  Suppose $\phi\in \mathcal{C}^2(\R_+)$ and $\psi\in \mathcal{C}^1(\R_+)$ are
subject to \eqref{da3} for some $\e>0$ and $k\geq 0$.  Now, consider the nonlinear wave equation with potential
\eqref{nl}.  Suppose the nonlinear term $F(u)$ satisfies \eqref{F} for some
\begin{equation}\label{pf}
p_n < p < 1 + \frac{2}{m} \:,
\end{equation}
where $p_n$ is the positive root of the quadratic \eqref{pn}.  Also, assume the potential term $V(r)$ is subject to
\eqref{V} for some $V_0>0$ and $\kappa>2$.  If $V_0,\e$ are sufficiently small, then \eqref{nl} admits a unique solution
$u\in \mathcal{C}^1(\Omega_T)$, where $T= +\infty$ in the supercritical case $k\geq 2/(p-1)$ and
\begin{equation}\label{lbf}
T \geq C\e^{-(p-1)/(2-k(p-1))}
\end{equation}
in the subcritical case $0\leq k< 2/(p-1)$.  Besides, the constant $C$ is independent of $\e$.
\end{theorem}

\begin{rem}
When it comes to initial data of \textit{subcritical} decay rate $0\leq k< 2/(p-1)$, the lower bound \eqref{lbf} for the
lifespan of solutions was obtained by Kubo \cite{Kb1}, still only for the special case $V(r)\equiv 0$ with $n$ odd.  Due
to a result of Takamura \cite{Ta2}, such a lower bound is known to be sharp when $V(r)\equiv 0$, regardless of the parity
of $n$.  As we shall prove later in this paper, it is actually sharp for any potential $V(r)$ that is merely non-positive
at infinity; see Theorem \ref{bu1}.
\end{rem}

\begin{rem}
When it comes to initial data of \textit{supercritical} decay rate $k\geq 2/(p-1)$, the existence of global solutions
persists in $n=3$ space dimensions as well.  In fact, a result of Strauss and Tsutaya \cite{ST} yields global
$\mathcal{C}^2$ solutions under similar assumptions that require more regularity, but not radial symmetry, of the initial
data and $V$. As we are going to show, however, our assumption \eqref{V} on the potential is not sufficient when $n=1,2$.
\end{rem}

To complement our existence result, Theorem \ref{et}, we shall also show that blow-up may occur for arbitrarily small
data under less favorable assumptions on either the initial data or the potential term.

In our first blow-up result, Theorem \ref{bu1}, blow-up occurs due to the slow decay rate of the initial data.  To merely
focus on the behavior of the initial data at infinity, we shall fix a constant $R>0$ and introduce the assumption
\begin{equation}\label{ba}
\phi(r)=0, \quad \psi(r)\geq \e r^{-k-1} \quad\quad \text{on\, $(R,\infty)$}
\end{equation}
for some $\e>0$ and $0\leq k< 2/(p-1)$.  For a potential $V(r)$ that is non-positive on $(R,\infty)$, we are then able to
establish the blow-up of solutions to \eqref{nl} when $F(u)= |u|^p$ or $|u|^{p-1} u$ for some $p>1$. Here, we also derive
an upper bound for the lifespan of local solutions which is similar to the lower bound given in \eqref{lbf}.

When it comes to initial data of noncompact support \eqref{ba}, there is a standard iteration method for proving blow-up
\cite{AT, As, Ta2}.  The underlying idea, which goes back to John \cite{Jo}, cannot be applied here directly, unless we
further restrict our initial data \eqref{ba} on the remaining interval $(0,R]$.  One way to get around this difficulty is
provided by Lemma \ref{cl}, a refinement of Keller's Comparison Theorem \cite{Ke} that would allow us to resort to the
standard iteration method.  Nevertheless, we use Lemma \ref{cl} to give a new and simpler method of proof which is based
on Glassey's ODE approach \cite{Gl2} for data of compact support.

In our second and last blow-up result, blow-up occurs due to the potential term.  Here, we remove our assumption of
radial symmetry and establish the following general

\begin{theorem}\label{bpg}
Let $n\geq 1$.  Suppose that $V\colon \R^n\to \R$ is continuous and that $-\Delta + V$ has a negative eigenvalue with a
positive eigenfunction which decays exponentially fast.  Let $\phi,\psi \geq 0$ be continuous functions of compact
support\footnote{We assume the data to be of compact support merely for the sake of simplicity.} and suppose that $\psi$
is not identically zero. If $(u, u_t) \in \mathcal{C}(H^1(\R^n) \times L^2(\R^n); [0,T))$ satisfies
\begin{equation}\label{pde}
\d_t^2 u - \Delta u + V(x)u = A|u|^p;\quad\quad u(x,0)= \phi(x),\quad \d_t u(x,0)= \psi(x)
\end{equation}
for some $A> 0$ and $p>1$, then $||u(\cdot \,,t)||_{L^q(\R^n)}$ blows-up in finite time for each $1\leq q\leq \infty$.
\end{theorem}

\begin{rem}
Similar blow-up results appear in \cite{ST, YZ2} but those require the potential term to be of one sign and also of rapid
decay at infinity.
\end{rem}

\begin{rem}
In section \ref{BP}, we give precise conditions on $V$ that ensure the applicability of this theorem.  Here, let us
merely remark that the eigenfunction corresponding to the first eigenvalue does have the desired properties under very
mild conditions on $V$.  In particular, the main hypothesis in this theorem is the presence of a negative eigenvalue.
This hypothesis holds for all potentials which behave like $-|x|^{-\kappa}$ at infinity for some $\kappa<2$.  Thus, the
decay assumption $\kappa>2$ in Theorem \ref{et} is almost necessary to ensure global solutions when $n\geq 3$.  When
$n=1,2$, the situation is slightly different because a negative eigenvalue may emerge even for potentials that are
rapidly decaying.
\end{rem}

\begin{rem}
If $V(x)\leq 0$ is a nonzero function of compact support, then $-\Delta + aV$ has a negative eigenvalue for all large
enough $a$.  Thus, one does need the potential term to be of small-amplitude in Theorem \ref{et}, as no sign condition is
imposed there.  When $n=3$, on the other hand, global solutions do exist for all non-negative potentials of compact
support \cite{GHK}.
\end{rem}

Finally, let us remark that our methods in this paper do not allow us to treat potentials which decay at the critical
rate $\kappa=2$. For that particular case, we refer the reader to \cite{PST1}.

The remaining of this paper is organized as follows.  Sections~\ref{pre} through \ref{bex} are devoted to the proof of
our existence result, Theorem \ref{et}.  In section~\ref{pre}, we review some facts about the homogeneous wave equation
and we introduce the weighted $L^\infty$ space in which solutions to \eqref{nl} are to be constructed. Section~\ref{slow}
contains certain estimates regarding our weight function which are needed in the proof of our existence result, while the
proof itself appears in section~\ref{bex}.  In section~\ref{BD}, we prove blow-up for initial data of subcritical decay,
while section~\ref{BP} settles our second blow-up result, Theorem \ref{bpg}.  Finally, section \ref{app} lists some facts
about the Riemann operator for the wave equation which were obtained in our previous work \cite{Ka1}.

\section{Preliminaries}\label{pre}
In this section, we prepare a few basic lemmas that will be needed in the proof of our existence theorem regarding the
nonlinear wave equation with potential
\begin{equation}\label{nlf}
\left\{
\begin{array}{rll}
\d_t^2 u -\d_r^2 u -\dfrac{n-1}{r} \cdot \d_r u = F(u) - V(r)\cdot u \quad\quad &\text{in\, $\Omega_T= \R_+\times
(0,T)$}\\
u(r,0) = \phi(r); \quad \d_t u(r,0)= \psi(r) \quad\quad &\text{in\, $\R_+$.}
\end{array}\right.
\end{equation}
Some of these lemmas depend on the parity of $n$ and, in particular, on the parameters
\begin{equation}\label{am}
(a,m) = \left\{
\begin{array}{ccl}
\left( 1\:,\: \frac{n-3}{2} \right) &&\text{if\, $n$ is odd} \\
\left( \frac{1}{2} \:,\: \frac{n-2}{2} \right) &&\text{if\, $n$ is even}
\end{array}\right.
\end{equation}
we shall frequently use in what follows.  We remark that $m\geq 1$ whenever $n\geq 4$, while the sum $a+m= (n-1)/2$ is
independent of the parity of $n$.

Recall that we seek a global solution to \eqref{nlf} for initial data of decay rate $k\geq 2/(p-1)$ and a local solution,
otherwise. There is no loss of generality in decreasing this decay rate as long as no lower bound on $k$ is contradicted.
In other words, we may take $k$ to be smaller than any quantity that exceeds $2/(p-1)$.  Now, our assumption \eqref{pf}
ensures that
\begin{equation*}
\frac{2}{p-1} < \frac{n-1}{2} \cdot p - 1 = (a+m)p -1,
\end{equation*}
as equality holds in the above inequality when $p= p_n$.  In particular, we may assume that
\begin{equation}\label{k1}
k < (a+m)p -1
\end{equation}
in what follows.  Similarly, one can readily check that
\begin{equation*}
p_n - 1 > \frac{4}{n+1} = \frac{2}{a+m+1}
\end{equation*}
and this allows us to additionally assume
\begin{equation}\label{k2}
k < \frac{n+1}{2} = a+m+1.
\end{equation}
Finally, it is convenient to decrease the decay rate $\kappa >2$ of the potential $V(r)$ so that
\begin{equation}\label{ka}
\kappa < m+2.
\end{equation}
We can do this without loss of generality when $m>0$, namely when $n\geq 4$.

Our plan is to construct a solution of \eqref{nlf} that is continuously differentiable and belongs to the Banach space
\begin{equation}\label{X}
X= \left\{ u(r,t)\in \mathcal{C}^1(\Omega_T) \::\: ||u|| < \infty \right\}, \quad\quad \Omega_T = \R_+\times (0,T).
\end{equation}
Here, the norm $||\cdot||$ is defined by
\begin{equation}\label{N}
||u|| = \sum_{j=0}^1 \sup_{(r,t)\in \Omega_T} \: |\d_r^j u(r,t)| \cdot r^{m-1+j} \br{r}^{1-j} \cdot W_k(r,t),
\end{equation}
where the weight function $W_k$ is of the form
\begin{equation}\label{W1}
W_k(r,t) = \br{t+r}^\mu \br{t-r}^\nu \left( 1 + \ln \frac{\br{t+r}}{\br{t-r}} \right)^{-\de_{k,m+a}}
\end{equation}
with $\mu= \min(k-m,a)$, $\nu= \max(k-m-a,0)$ and $\de_{k,m+a}$ the usual Kronecker delta.  This weighted norm is partly
dictated by our previous work \cite{Ka1} on the homogeneous problem
\begin{equation}\label{he}
\left\{
\begin{array}{rll}
\d_t^2 u_0 - \d_r^2 u_0  - \dfrac{n-1}{r} \cdot \d_r u_0 = 0 \quad\quad &\text{in\: $\R_+^2= (0,\infty)^2$} \\
u_0(r,0) = \phi(r);\:\: \d_t u_0(r,0) = \psi(r) &\text{in\: $\R_+$.}
\end{array}\right.
\end{equation}

\begin{lemma}\label{dec}
Let $n\geq 4$ be an integer and define $a,m$ by \eqref{am}.  Suppose that $\phi\in \mathcal{C}^2(\R_+)$ and $\psi\in
\mathcal{C}^1(\R_+)$ are subject to \eqref{da3} for some $\e>0$ and some $0\leq k <(n+1)/2$. Then, the homogeneous
equation \eqref{he} admits a unique solution $u_0\in \mathcal{C}^1(\R_+^2)$ which satisfies
\begin{equation}\label{en}
|\d_r^j u_0| + |\d_t^j u_0| \leq C_0(k,n) \cdot \e r^{1-m-j} \cdot \br{t-r}^{-j} \,\br{t+r}^{j-1} \cdot W_k(r,t)^{-1}
\end{equation}
when $j=0,1$.  In particular, $u_0$ is in the Banach space \eqref{X} and we have $||u_0||\leq C_0\e$.
\end{lemma}

\begin{proof}
Decay estimates for the solution to the homogeneous wave equation \eqref{he} appear in Theorem 1.1 of \cite{Ka1}.
Although no restrictions were imposed there on the decay rate $k$ of the initial data, we shall only need to treat decay
rates $0\leq k< (n+1)/2$ here; see \eqref{k2}.  Under our assumption that $n\geq 4$, such decay rates fall in the range
$0\leq k< n-1$.  According to Theorem 1.1 in \cite{Ka1} then, \eqref{he} has a unique solution $u_0\in
\mathcal{C}^1(\R_+^2)$ which satisfies \eqref{en}.  This also implies
\begin{equation*}
|\d_r^j u_0| \leq C_0(k,n) \cdot \e r^{1-m-j} \,\br{r}^{j-1} \cdot W_k(r,t)^{-1}
\end{equation*}
when $j=0,1$, so the definition \eqref{N} of our norm allows us to deduce that $||u_0||\leq C_0\e$.
\end{proof}

The main purpose of our previous work \cite{Ka1} was to study the Riemann operator $L$ for the wave equation in the
radial case.  Section \ref{app} lists some of the estimates we established there, as those are also useful in treating
the nonlinear wave equation \eqref{nlf}.  In fact, the standard Duhamel principle allows us to obtain the following

\begin{lemma}\label{new}
Let $L$ denote the Riemann operator of Lemma \ref{hs}.  For a function $G$ of two variables, we define the Duhamel
operator $\mathscr{L}$ as
\begin{equation}\label{Du}
[\mathscr{L}G](r,t) = \int_0^t [LG(\cdot\,,\tau)](r,t-\tau) \:d\tau.
\end{equation}
When $G\in \mathcal{C}^1(\Omega_T)$, one then has $\mathscr{L}G\in \mathcal{C}^1(\Omega_T)$ and this function provides a
solution to
\begin{equation*}
\left( \d_t^2 - \d_r^2 - \frac{n-1}{r} \cdot \d_r \right) [\mathscr{L}G](r,t) = G(r,t) \quad\quad \text{in\, $\Omega_T=
\R_+\times (0,T)$}
\end{equation*}
when zero initial data are imposed.
\end{lemma}

\begin{proof}
Our assertions follow easily by means of Lemma \ref{hs} and a simple computation.
\end{proof}

\begin{prop}\label{Les}
Let $n\geq 4$ be an integer and define $a,m$ by \eqref{am}.  Suppose $G\in \mathcal{C}^1(\Omega_T)$ satisfies the
singularity condition
\begin{equation}\label{sg2}
G(\la,\tau) = O \left( \la^{-2m-2+\de} \right) \quad \text{as $\la \rightarrow 0$}
\end{equation}
for some fixed $\de>0$.  With $D= (\d_r, \d_t)$ and $\la_\pm = t-\tau\pm r$, one then has
\begin{align*}
|D^\b [\mathscr{L}G](r,t)|
&\leq Cr^{j-|\b|-m-a} \int_0^t \int_{|\la_-|}^{\la_+} \frac{\la^{m-j+1}}{(\la-\la_-)^{1-a}} \cdot \sum_{s=0}^j
\la^s \,|\d_\la^s G(\la,\tau)| \:d\la\,d\tau \\
&\quad + Cr^{j-|\b|-m} \int_0^{t-r} \int_0^{\la_-} \frac{\la^{2m-j+1}}{\la_-^m \la_+^a \,(\la_- -\la)^{1-a}} \cdot
\sum_{s=0}^j \la^s \,|\d_\la^s G(\la,\tau)| \:d\la\,d\tau\\
&\quad + Cr^{j-|\b|-m-a} \int_{\max(t-2r,0)}^t  |\la_\pm|^{a+m-j+1} \left[ \sum_{s=0}^{j-1} \la^s \,|\d_\la^s G(\la,\tau)|
\right]_{\la= |\la_\pm|} d\tau
\end{align*}
whenever $|\b|\leq j\leq 1$.  Besides, the constant $C$ is independent of $r,t$.
\end{prop}

\begin{proof}
Since the integrand in \eqref{Du} vanishes when $\tau= t$, a direct differentiation gives
\begin{align}\label{mm}
D^\b [\mathscr{L}G]
&= \int_{\max(t-2r,0)}^t D^\b [LG(\cdot\,,\tau)](r,t-\tau) \:d\tau + \int_0^{\max(t-2r,0)} D^\b [LG(\cdot\,,\tau)]
(r,t-\tau) \:d\tau \notag\\
&\equiv A + B.
\end{align}
To treat the first integral $A$, we note that $t-\tau\leq 2r$ within the region of integration.  This allows us to invoke
Lemma \ref{Lex} to find that
\begin{align*}
A &\leq C_1(n) \cdot r^{-m-a} \int_{\max(t-2r,0)}^t \int_{|\la_-|}^{\la_+} \frac{\la^{m-|\b|+1}}{(\la-\la_-)^{1-a}}
\cdot \sum_{s=0}^j \la^s \,|\d_\la^s G(\la,\tau)| \:d\la\,d\tau \\
&\quad + C_2(n) \cdot r^{-m-a} \int_{\max(t-2r,0)}^{t-r} \int_0^{\la_-} \frac{\la^{2m+1}}{\la_-^{m+|\b|}
\,(\la_- -\la)^{1-a}} \cdot \sum_{s=0}^j \la^s \,|\d_\la^s G(\la,\tau)| \:d\la\,d\tau\\
&\quad + C_1(n) \cdot r^{-m-a} \int_{\max(t-2r,0)}^t  |\la_\pm|^{a+m-|\b|+1} \left[ \sum_{s=0}^{j-1} \la^s
\,|\d_\la^s G(\la,\tau)| \right]_{\la= |\la_\pm|} d\tau \\
&\equiv A_1 + A_2 + A_3
\end{align*}
with $\la_\pm = t-\tau\pm r$.  Let us merely concern ourselves with the middle term $A_2$, as the other terms are easier
to handle. Since $\la_+\leq 3r$ within the region of integration, we get
\begin{equation*}
A_2 \leq Cr^{-m} \int_{\max(t-2r,0)}^{t-r} \int_0^{\la_-} \frac{\la^{2m-|\b|+1}}{\la_-^m \la_+^a \,(\la_- -\la)^{1-a}}
\cdot \sum_{s=0}^j \la^s \,|\d_\la^s G(\la,\tau)| \:d\la\,d\tau.
\end{equation*}
Moreover, one has $\la\leq \la_+\leq 3r$ whenever $t-\tau\leq 2r$, so our assumption $|\b|\leq j$ implies
\begin{equation*}
A_2 \leq Cr^{j-|\b|-m} \int_0^{t-r} \int_0^{\la_-} \frac{\la^{2m-j+1}}{\la_-^m \la_+^a \,(\la_- -\la)^{1-a}} \cdot
\sum_{s=0}^j \la^s \,|\d_\la^s G(\la,\tau)| \:d\la\,d\tau.
\end{equation*}
In particular, the desired estimate is satisfied by the first integral in \eqref{mm}.

Returning to \eqref{mm}, we now focus on the second integral $B$.  Since $t-\tau\geq 2r$ within the region of
integration, we may apply Lemma \ref{Lin} to similarly get
\begin{align*}
B &\leq C_1'(n) \cdot r^{j-|\b|-m-a} \int_0^{\max(t-2r,0)} \int_{\la_-}^{\la_+} \frac{\la^{m-j+1}}{(\la-\la_-)^{1-a}}
\cdot \sum_{s=0}^j \la^s \,|\d_\la^s G(\la,\tau)| \:d\la\,d\tau \\
&\quad + C_2'(n) \cdot r^{j-|\b|-m} \int_0^{\max(t-2r,0)} \int_0^{\la_-} \frac{\la^{2m+1}}{\la_-^{j+m+a}
\,(\la_- -\la)^{1-a}} \cdot \sum_{s=0}^j \la^s \,|\d_\la^s G(\la,\tau)| \:d\la\,d\tau\\
&\equiv B_1 + B_2.
\end{align*}
It is clear that the first integral satisfies the desired estimate, so we need only worry about the second. Since $\la_-=
t-\tau-r$ and $\la_+= t-\tau+r$ are equivalent when $t-\tau\geq 2r$, we find
\begin{equation*}
B_2 \leq Cr^{j-|\b|-m} \int_0^{t-r} \int_0^{\la_-} \frac{\la^{2m-j+1}}{\la_-^m \la_+^a \,(\la_- -\la)^{1-a}} \cdot
\sum_{s=0}^j \la^s \,|\d_\la^s G(\la,\tau)| \:d\la\,d\tau.
\end{equation*}
In particular, the desired estimate is satisfied by $B_2$ as well and the proof is complete.
\end{proof}

\begin{lemma}\label{bas}
Suppose $u$ belongs to the Banach space \eqref{X} and let $p>1$.  Assuming \eqref{F} and \eqref{V}, one then has
\begin{equation}\label{ma1}
\sum_{s=0}^{j_0} \la^s \:|\d_\la^s F(u(\la,\tau))| \leq 2Ap ||u||^p \cdot\la^{j-mp} \br{\la}^{j_0-j} \cdot
W_k(\la,\tau)^{-p}
\end{equation}
and also
\begin{equation}\label{ma2}
\sum_{s=0}^j \la^s \:|\d_\la^s (V(\la)\cdot u(\la,\tau))| \leq 4V_0 ||u|| \cdot \la^{j-1-m} \br{\la}^{1-\kappa} \cdot
W_k(\la,\tau)^{-1}
\end{equation}
whenever $0\leq j, j_0\leq 1$ and $(\la,\tau)\in \Omega_T$.
\end{lemma}

\begin{proof}
Because of our assumption \eqref{F}, the fundamental theorem of calculus ensures that
\begin{equation*}
|F(u)| \leq A\cdot |u|^p, \quad\quad |F'(u)|\leq Ap\cdot |u|^{p-1}.
\end{equation*}
In particular, it ensures that
\begin{equation*}
|\d_\la^s F(u(\la,\tau))| \leq Ap^s \cdot |u(\la,\tau)|^{p-s} \cdot |\d_\la u(\la,\tau)|^s, \quad\quad s=0,1.
\end{equation*}
Recalling the definition \eqref{N} of our norm, the last equation easily leads to
\begin{equation*}
\la^s \,|\d_\la^s F(u(\la,\tau))|\leq Ap^s ||u||^p \cdot\la^{p-mp}\br{\la}^{s-p}\cdot W_k(\la,\tau)^{-p}, \quad\quad
s=0,1.
\end{equation*}
In view of our assumptions that $j\leq 1< p$, this also implies
\begin{equation*}
\sum_{s=0}^{j_0} \la^s \,|\d_\la^s F(u(\la,\tau))| \leq Ap ||u||^p \cdot \sum_{s=0}^{j_0} \la^{j-mp} \br{\la}^{s-j} \cdot
W_k(\la,\tau)^{-p}, \quad\quad j_0=0,1.
\end{equation*}
Moreover, $s-j\leq j_0-j$ within the last sum, so our first assertion \eqref{ma1} follows.

Since our second assertion \eqref{ma2} is easier to establish, we shall omit the details.
\end{proof}

\section{A Priori Estimates}\label{slow}
Our main goal in this section is to establish the following

\begin{theorem}\label{bes}
Let $n\geq 4$ be an integer and define $a,m$ by \eqref{am}.  Suppose $F(u)$ and $V(r)$ are subject to \eqref{F} and
\eqref{V}, respectively.  Suppose the decay rates $k,\kappa$ are subject to \eqref{k1} through \eqref{ka} and that the
condition
\begin{equation}\label{con}
2 - mp+ m > 0
\end{equation}
holds.  Define the norm $||\cdot||$ as in \eqref{N} and the function $\Phi_k$ by the formula
\begin{equation}\label{Phi}
\Phi_k(t+r) = \br{t+r}^{\max(2-k(p-1),0)}.
\end{equation}
Let $\mathscr{L}$ be the Duhamel operator \eqref{Du} and assume $u$ belongs to the Banach space \eqref{X}.  With $D=
(\d_r, \d_t)$, one then has
\begin{equation}\label{mb1}
|D^\b [\mathscr{L} F(u)](r,t)| \leq C||u||^p \cdot r^{1-|\b|-m} \br{r}^{|\b|-1} \cdot W_k(r,t)^{-1} \cdot \Phi_k(t+r)
\end{equation}
and also
\begin{equation}\label{mb2}
|D^\b [\mathscr{L} (Vu)](r,t)| \leq CV_0 ||u||\cdot r^{1-|\b|-m} \br{r}^{|\b|-1} \cdot W_k(r,t)^{-1}
\end{equation}
as long as $|\b|\leq 1$ and $(r,t)\in \Omega_T$.  Besides, the constant $C$ is independent of $r, t$.
\end{theorem}

Before we turn to the proof of this theorem, we first need to study certain integrals which will arise in the course of
the proof.  Those involve our weight function \eqref{W1} and some other parameters we have introduced \eqref{am}.
Throughout this section, in particular, we assume
\begin{equation}\label{ks}
0\leq k < m+a+1;\quad a>0; \quad \mu= \min(k-m,a);\quad \nu= \max(k-m-a,0).
\end{equation}
For future reference, let us remark that $\mu+ \nu= k-m$ and that $0\leq \nu < 1$.

\begin{lemma}\label{fa1}
Let $b, y\geq 0$ be arbitrary.  When $\kappa>2$ and $\nu< 1$, one then has
\begin{equation}\label{fa1a}
I_{1b} \equiv \int_0^y \left( 1+ \ln \frac{\br{y}}{\br{x}} \right)^b dx \leq C(b) \cdot y
\end{equation}
as well as
\begin{equation}\label{fa1b}
I_2 \equiv \int_{-y}^y \br{x+y}^{1-\kappa} \cdot \br{x}^{-\nu} \left( 1+ \ln \frac{\br{y}}{\br{x}} \right)^b dx \leq
C(\kappa, \nu, b) \cdot \br{y}^{-\nu}.
\end{equation}
\end{lemma}

\begin{proof}
The given integrals are increasing in $b$, so we may assume that $b$ is an integer.  Let us first focus on
\eqref{fa1a}. Since $I_{10} = y$ and since an integration by parts gives
\begin{equation*}
I_{1b} = y + \int_0^y \frac{bx}{1+x} \cdot \left( 1+ \ln \frac{\br{y}}{\br{x}} \right)^{b-1} dx \leq y + b I_{1,b-1},
\end{equation*}
the validity of \eqref{fa1a} follows by induction.

Next, we turn to \eqref{fa1b}.  Here, $x+y$ is equivalent to $y$ for each $-y/2 \leq x \leq y$ and $\br{x}$ is equivalent
to $\br{y}$ for the remaining part $-y\leq x\leq -y/2$, so we get
\begin{equation*}
I_2 \leq C\br{y}^{1-\kappa} \int_{-y/2}^y \br{x}^{-\nu} \left( 1+ \ln \frac{\br{y}}{\br{x}} \right)^b dx + C\br{y}^{-\nu}
\int_{-y}^{-y/2} \br{x+y}^{1-\kappa} \:dx.
\end{equation*}
Recalling our assumption that $\kappa>2$, we then arrive at
\begin{equation*}
I_2 \leq C\br{y}^{1-\kappa+\kappa-2} \int_{-y/2}^y \br{x}^{-\nu} dx + C\br{y}^{-\nu}.
\end{equation*}
Since we also have $\nu< 1$ by assumption, the desired estimate \eqref{fa1b} follows trivially.
\end{proof}

\begin{lemma}\label{fa2}
Let $y\geq 0$ and $p>1$.  Assuming \eqref{k1}, \eqref{con} and \eqref{ks}, one has
\begin{equation}\label{fa2a}
J_1 \equiv \int_{-y}^y (x+y)^{1-mp+m} \cdot \br{x}^{-\nu p} \left( 1+ \ln \frac{\br{y}}{\br{x}} \right)^{p\de_{k,m+a}} dx
\leq C\Phi_k(y) \cdot \br{y}^{\mu p +m-k},
\end{equation}
where $\Phi_k(y) = \br{y}^{\max(2-k(p-1),0)}$ is given by \eqref{Phi} and $C$ is independent of $y$.
\end{lemma}

\begin{proof}
If it happens that $0\leq y\leq 1$, then we easily find
\begin{equation*}
J_1 \leq C\int_{-y}^y (x+y)^{1-mp+m} \:dx \leq Cy^{2-mp+m} \leq C
\end{equation*}
because $2-mp+m>0$ by \eqref{con}.  This does imply the desired \eqref{fa2a} when $y$ is bounded.

Assume now that $y\geq 1$. Since $x+y$ is equivalent to $\br{y}$ for each $-y/2 \leq x \leq y$ and $\br{x}$ is equivalent
to $\br{y}$ for the remaining part $-y\leq x\leq -y/2$, we get
\begin{equation*}
J_1 \leq C\br{y}^{1-mp+m} \int_{-y/2}^y \br{x}^{-\nu p} \left( 1+ \ln \frac{\br{y}}{\br{x}} \right)^{p\de_{k,m+a}} dx +
C\br{y}^{-\nu p} \int_{-y}^{-y/2} (x+y)^{1-mp+m} \:dx.
\end{equation*}
Moreover, $2-mp+m$ is positive by \eqref{con}, so this actually gives
\begin{equation}\label{J11}
J_1 \leq C\br{y}^{1-mp+m} \int_{-y/2}^y \br{x}^{-\nu p} \left( 1+ \ln \frac{\br{y}}{\br{x}} \right)^{p\de_{k,m+a}} dx +
C\br{y}^{2-mp+m -\nu p}.
\end{equation}
Recalling that $\mu= \min(k-m,a)$ and $\nu= \max(k-m-a,0)$, we shall consider two cases.

\Case{1} When $k\leq m+a$, we have $\mu = k-m$ and $\nu= 0$, so the last equation reads
\begin{equation*}
J_1 \leq C\br{y}^{1-mp+m} \int_{-y/2}^y \left( 1+ \ln \frac{\br{y}}{\br{x}} \right)^{p\de_{k,m+a}} dx +C\br{y}^{2-mp+m}.
\end{equation*}
Once we now employ \eqref{fa1a} to treat the integral, we find
\begin{equation*}
J_1 \leq C\br{y}^{2-mp+m} = C\br{y}^{2-k(p-1)} \cdot \br{y}^{\mu p +m-k}
\end{equation*}
since $\mu= k-m$ for this case. In view of the definition of $\Phi_k$, the desired \eqref{fa2a} follows.

\Case{2} When $k>m+a$, we have $\mu =a$ and $\nu = k-m-a$, while equation \eqref{J11} reads
\begin{equation*}
J_1 \leq C\br{y}^{1-mp+m} \int_{-y/2}^y \br{x}^{-\nu p} dx + C\br{y}^{2-mp+m -\nu p}.
\end{equation*}
Moreover, $\nu= k-m-a$ for this case, so we also have the identity
\begin{equation*}
- \nu p = [2-k(p-1)] + (m+a)p - k-2.
\end{equation*}
In view of the definition of $\Phi_k$, the last two equations combine to give
\begin{equation*}
J_1 \leq  C\Phi_k(y)\cdot \br{y}^{1-mp+m} \int_{-y/2}^y \br{x}^{(m+a)p-k-2} \:dx + C\Phi_k(y)\cdot \br{y}^{ap+m-k}.
\end{equation*}
Since $(m+a)p-k-1>0$ by our assumption \eqref{k1}, we thus obtain the estimate
\begin{equation*}
J_1 \leq C\Phi_k(y) \cdot \br{y}^{ap+m-k}.
\end{equation*}
In particular, we obtain the desired estimate \eqref{fa2a} because $\mu= a$ for this case.
\end{proof}

The proof of the following fact is almost identical to that of Lemma 3.6 in \cite{Ts1}, so we are going to omit it.

\begin{lemma}\label{ts1}
Let $(r,t)\in \R_+^2$ be arbitrary.  Assuming that $a>0$, one has
\begin{equation*}
\int_{|t-r|}^{t+r} \frac{\br{\la}^b d\la}{(r-t+\la)^{1-a}} \leq \left\{
\begin{array}{llcc}
Cr^a \br{t+r}^b &&\text{if} & b > -a \\
Cr^a \br{t+r}^{-a} \left( 1+ \ln \dfrac{\br{t+r}}{\br{t-r}} \right) &&\text{if} & b = -a \\
Cr^a \br{t+r}^{-a} \br{t-r}^{a+b} &&\text{if} & b < -a \\
\end{array}\right\}
\end{equation*}
for some constant $C$ depending solely on $a$ and $b$.
\end{lemma}

\begin{lemma}\label{I1}
Let $W_k$ be the weight function \eqref{W1} and $\kappa>2$.  Assuming \eqref{ks}, one has
\begin{equation*}
\mathcal{I}_1 \equiv \int_0^t \int_{|\la_-|}^{\la_+} \frac{\br{\la}^{1-\kappa} \cdot
W_k(\la,\tau)^{-1}}{(\la-\la_-)^{1-a}} \:\:d\la\,d\tau \leq Cr^a W_k(r,t)^{-1},
\end{equation*}
where $\la_\pm = t-\tau\pm r$ and the constant $C$ is independent of $r, t$.
\end{lemma}

\begin{proof}
Let us recall the definition \eqref{W1} of our weight function $W_k$ and write
\begin{align*}
\mathcal{I}_1 &= \int_0^t \int_{|t-\tau-r|}^{t-\tau+r} \frac{\br{\la}^{1-\kappa} \cdot \br{\la+\tau}^{-\mu}}{(r-t+ \tau+
\la)^{1-a}} \cdot \br{\la-\tau}^{-\nu} \left( 1+ \ln \frac{\br{\la+\tau}}{\br{\la-\tau}} \right)^{\de_{k,m+a}}
d\la\,d\tau.
\end{align*}
Changing variables by $x= \la - \tau$ and $y = \la+\tau$, we then get
\begin{equation*}
\mathcal{I}_1 \leq C\int_{|t-r|}^{t+r} \frac{\br{y}^{-\mu}}{(r-t+y)^{1-a}} \int_{-y}^y \br{x+y}^{1-\kappa} \cdot
\br{x}^{-\nu} \left( 1+ \ln \frac{\br{y}}{\br{x}} \right)^{\de_{k,m+a}} dx\,dy.
\end{equation*}
Once we now employ Lemma \ref{fa1} to treat the inner integral, we arrive at
\begin{equation*}
\mathcal{I}_1 \leq C\int_{|t-r|}^{t+r} \frac{\br{y}^{-\mu-\nu}}{(r-t+y)^{1-a}} \:\:dy = C\int_{|t-r|}^{t+r}
\frac{\br{y}^{m-k}}{(r-t+y)^{1-a}} \:\:dy
\end{equation*}
because $\mu+ \nu= k-m$ by \eqref{ks}.  To finish the proof, it thus suffices to show
\begin{equation*}
\int_{|t-r|}^{t+r} \frac{\br{y}^{m-k}}{(r-t+y)^{1-a}} \:\:dy \leq Cr^a W_k(r,t)^{-1}.
\end{equation*}
In other words, it suffices to show
\begin{equation*}
\int_{|t-r|}^{t+r} \frac{\br{y}^{m-k}}{(r-t+y)^{1-a}} \:\:dy \leq Cr^a\br{t+r}^{-\mu} \br{t-r}^{-\nu} \left( 1+ \ln
\frac{\br{t+r}}{\br{t-r}} \right)^{\de_{k,m+a}}
\end{equation*}
with $\mu= \min(k-m,a)$, $\nu= \max(k-m-a,0)$ and $\de_{k,m+a}$ the usual Kronecker delta.  Since this is precisely the
estimate provided by the previous lemma, the proof is complete.
\end{proof}

Repeating the above proof but using Lemma \ref{fa2} to treat the inner integral, one obtains

\begin{lemma}\label{J1}
Let $W_k$ be as in \eqref{W1} and $p>1$.  Assuming \eqref{k1}, \eqref{con} and \eqref{ks}, one has
\begin{equation*}
\mathcal{J}_1 \equiv \int_0^t \int_{|\la_-|}^{\la_+} \frac{\la^{1-mp+m} \cdot W_k(\la,\tau)^{-p}}{(\la-\la_-)^{1-a}}
\:\:d\la\,d\tau \leq Cr^a W_k(r,t)^{-1} \cdot \Phi_k(t+r),
\end{equation*}
where $\la_\pm = t-\tau\pm r$, $\Phi_k$ is given by \eqref{Phi} and the constant $C$ is independent of $r, t$.
\end{lemma}

The proof of the following fact is almost identical to that of Lemma 3.7 in \cite{Ts1}, so we are going to omit it.

\begin{lemma}\label{ts2}
Given constants $a>0$, $b\geq 0$ and $c> -b-1$, one has
\begin{equation*}
\int_0^{t-r} \frac{\la^b \br{\la}^c \:d\la}{(t-r-\la)^{1-a}} \leq C(t-r)^{a+b} \cdot \br{t-r}^c
\end{equation*}
whenever $t\geq r> 0$.  Besides, the constant $C$ depends solely on $a$, $b$ and $c$.
\end{lemma}

\begin{lemma}\label{I2}
Define $W_k$ by \eqref{W1} and let $\kappa>2$.  Assuming \eqref{ka} and \eqref{ks} with $m\geq a$, one then has
\begin{equation*}
\mathcal{I}_2 \equiv \int_0^{t-r} \int_0^{\la_-} \frac{\la^m \br{\la}^{1-\kappa} \cdot W_k(\la,\tau)^{-1}}{\la_-^m
\la_+^a \,(\la_- -\la)^{1-a}} \:\:d\la\,d\tau \leq CW_k(r,t)^{-1},
\end{equation*}
where $\la_\pm = t-\tau \pm r$ and the constant $C$ is independent of $r,t$.
\end{lemma}

\begin{proof}
Here, the factor $\la_+^a$ in the denominator has to be treated carefully, so we shall need to divide our analysis
into two cases. Before we do this, however, let us first note that
\begin{equation}\label{fa}
\frac{\la}{\:\:\la_-} = \frac{\la}{t-\tau-r}  \leq \frac{C(\la + \tau)}{t-r} \:.
\end{equation}
This holds if $\tau \geq (t-r)/2$ and $0\leq \la\leq t-\tau-r$, in which case $\la+\tau$ is equivalent to $t-r$, but it
also holds if $0\leq \tau \leq (t-r)/2$, in which case $t-r-\tau$ is equivalent to $t-r$.

\Case{1} When $t\leq 2r$, the fact that $\la_+= t-\tau+r\geq r$ combines with \eqref{fa} to yield
\begin{equation*}
\mathcal{I}_2 \leq \frac{C}{r^a(t-r)^m} \int_0^{t-r} \int_0^{t-\tau -r} \frac{\br{\la}^{1-\kappa} \cdot (\la+\tau)^m
\cdot W_k(\la,\tau)^{-1}}{(t-\tau - r-\la)^{1-a}} \:\:d\la\,d\tau.
\end{equation*}
Recalling our definition \eqref{W1}, let us change variables by $x= \la-\tau$ and $y= \la+\tau$ to write
\begin{equation*}
\mathcal{I}_2 \leq \frac{C}{r^a(t-r)^m} \int_0^{t-r} \frac{y^m \br{y}^{-\mu}}{(t-r-y)^{1-a}}
\int_{-y}^y\br{x+y}^{1-\kappa} \cdot \br{x}^{-\nu} \left( 1+ \ln \dfrac{\br{y}}{\br{x}} \right)^{\de_{k,m+a}} dx\,dy.
\end{equation*}
Once we now employ Lemma \ref{fa1} to treat the inner integral, we find
\begin{equation}\label{In2}
\mathcal{I}_2 \leq \frac{C}{r^a(t-r)^m} \int_0^{t-r} \frac{y^m \br{y}^{-\mu-\nu}}{(t-r-y)^{1-a}} \:\:dy.
\end{equation}
Since $\mu\leq a$ and $\nu< 1$ by \eqref{ks}, our assumption $a\leq m$ in this theorem gives
\begin{equation*}
m- \mu -\nu \geq a-\mu -\nu \geq -\nu > -1.
\end{equation*}
As long as the last inequality holds, however, Lemma \ref{ts2} provides the estimate
\begin{equation}\label{In3}
\int_0^{t-r} \frac{y^m \br{y}^{-\mu-\nu}}{(t-r-y)^{1-a}} \:\:dy \leq C(t-r)^{m+a} \br{t-r}^{-\mu-\nu}.
\end{equation}
Combining this with \eqref{In2}, we may deduce the desired estimate once we show that
\begin{equation}\label{ne}
r^{-a} (t-r)^a \br{t-r}^{-\mu-\nu} \leq C\br{t+r}^{-\mu} \br{t-r}^{-\nu} \quad\quad\text{when $t\leq 2r$.}
\end{equation}
If $t\leq 2r$ and $r\leq 1$, this is easy to see because $t-r\leq r$ and $\br{t-r}$ is equivalent to $\br{t+r}$.  If
$t\leq 2r$ and $r\geq 1$, on the other hand, $r$ is equivalent to $\br{t+r}$ and we similarly get
\begin{equation*}
r^{-a} (t-r)^a \br{t-r}^{-\mu-\nu} \leq C\br{t+r}^{-a} \br{t-r}^{a-\mu-\nu} \leq C\br{t+r}^{-\mu} \br{t-r}^{-\nu}
\end{equation*}
because $\mu\leq a$ by \eqref{ks}.

\Case{2} When $t\geq 2r$, it is convenient to express the given integral as the sum of
\begin{equation*}
\mathcal{I}_{21} = \int_0^{2(t-r)/3} \int_0^{\la_-}  \frac{\la^m \br{\la}^{1-\kappa} \cdot W_k(\la, \tau)^{-1}}{\la_-^m
\la_+^a \,(\la_- -\la)^{1-a}} \:\:d\la\,d\tau
\end{equation*}
and
\begin{equation*}
\mathcal{I}_{22} = \int_{2(t-r)/3}^{t-r} \int_0^{\la_-}  \frac{\la^m \br{\la}^{1-\kappa} \cdot W_k(\la,
\tau)^{-1}}{\la_-^m \la_+^a \,(\la_- -\la)^{1-a}} \:\:d\la\,d\tau.
\end{equation*}
To treat $\mathcal{I}_{21}$, we proceed as in Case 1 with the inequality $\la_+= t-\tau+r\geq (t-r)/3$ instead of our
previous $\la_+ \geq r$.  Analogously to \eqref{In2}, we now establish
\begin{equation*}
\mathcal{I}_{21} \leq  \frac{C}{(t-r)^{m+a}} \int_0^{t-r} \frac{y^m \br{y}^{-\mu-\nu}}{(t-r-y)^{1-a}} \:\:dy \leq
C\br{t-r}^{-\mu-\nu}
\end{equation*}
because of \eqref{In3}.  This does imply the desired estimate whenever $t\geq 2r$.

To treat $\mathcal{I}_{22}$, we use the inequality $\la_+ \geq \la_-$ to first obtain
\begin{equation*}
\mathcal{I}_{22} \leq \int_{2(t-r)/3}^{t-r} \int_0^{\la_-} \frac{\la^m \br{\la}^{1-\kappa} \cdot
W_k(\la,\tau)^{-1}}{\la_-^{m+a} \,(\la_- -\la)^{1-a}} \:\:d\la\,d\tau.
\end{equation*}
Here, each of $\tau \pm \la$ is equivalent to $\tau$ within the region of integration because
\begin{equation*}
2\la \leq 2\la_- = 2(t-r-\tau) \leq \tau
\end{equation*}
whenever $\tau \geq 2(t-r)/3$.  In particular, each of $\tau \pm \la$ is equivalent to $t-r$ and so
\begin{equation*}
W_k(\la,\tau) = \br{\tau+\la}^\mu \br{\tau-\la}^\nu \left( 1+ \ln \frac{\br{\tau+\la}}{\br{\tau-\la}}
\right)^{-\de_{k,m+a}}
\end{equation*}
is equivalent to $\br{t-r}^{\mu+\nu}$.  Keeping this in mind, we then trivially get
\begin{equation*}
\mathcal{I}_{22} \leq C\br{t-r}^{-\mu-\nu} \int_0^{t-r} \int_0^{\la_-} \frac{\la^m \br{\la}^{1-\kappa}}{\la_-^{m+a}
\,(\la_- -\la)^{1-a}} \:\:d\la\,d\tau.
\end{equation*}
Since our next lemma shows the last integral is finite, the proof is finally complete.
\end{proof}

\begin{lemma}\label{I2s}
Let $m\geq a>0$ and $\kappa>2$.  Assuming \eqref{ka}, one has
\begin{equation*}
\mathbb{I} \equiv \int_0^{t-r} \int_0^{\la_-} \frac{\la^m \br{\la}^{1-\kappa}}{\la_-^{m+a} \,(\la_- -\la)^{1-a}}
\:\:d\la\,d\tau \leq C(a,m,\kappa).
\end{equation*}
\end{lemma}

\begin{proof}
Let us change variables by $y= \la_- = t-r-\tau$ and write
\begin{equation*}
\mathbb{I} = \int_0^{t-r}\int_{y/2}^y \frac{\la^m \br{\la}^{1-\kappa} \:d\la\,dy}{y^{m+a} \,(y-\la)^{1-a}} +
\int_0^{t-r}\int_0^{y/2} \frac{\la^m \br{\la}^{1-\kappa} \:d\la\,dy}{y^{m+a} \,(y-\la)^{1-a}} \equiv \mathbb{I}_1 +
\mathbb{I}_2.
\end{equation*}
For the former integral, the equivalence of $\la$ with $y$ easily leads to
\begin{equation*}
\mathbb{I}_1\leq C\int_0^{t-r} y^{-a}\br{y}^{1-\kappa}\int_{y/2}^y (y-\la)^{a-1} \:d\la\,dy = C\int_0^{t-r}
\br{y}^{1-\kappa}\:dy\leq C
\end{equation*}
because $a>0$ and $\kappa>2$.  For the latter integral, the equivalence of $y-\la$ with $y$ gives
\begin{equation*}
\mathbb{I}_2 \leq C\int_0^{t-r} y^{-m-1} \int_0^{y/2} \la^m \br{\la}^{1-\kappa} \:d\la\,dy.
\end{equation*}
Once we now consider the regions $y\geq 1$ and $y\leq 1$ separately, we find
\begin{equation*}
\mathbb{I}_2 \leq C\int_{\min(t-r,1)}^{t-r} \br{y}^{-m-1} \int_0^{y/2} \br{\la}^{m+1-\kappa} \:d\la\,dy +
C\int_0^{\min(t-r,1)} y^{-m-1} \int_0^{y/2} \la^m \:d\la\,dy
\end{equation*}
because $m>0$.  Using our assumption \eqref{ka} that $m+2-\kappa>0$, we then arrive at
\begin{equation*}
\mathbb{I}_2 \leq C\int_{\min(t-r,1)}^{t-r} \br{y}^{1-\kappa} \:dy + C\int_0^{\min(t-r,1)} \:dy.
\end{equation*}
Since we also have $\kappa>2$, the integrals on the right hand side are bounded, indeed.
\end{proof}

\begin{lemma}\label{J2}
Let $W_k$ be as in \eqref{W1}.  Assume \eqref{k1}, \eqref{con} and \eqref{ks} with $m\geq a$.  In the case that $p>1$,
one then has
\begin{equation*}
\mathcal{J}_2 \equiv \int_0^{t-r} \int_0^{\la_-} \frac{\la^{1-mp+2m} \cdot W_k(\la,\tau)^{-p}}{\la_-^m \la_+^a \,(\la_-
-\la)^{1-a}} \:\:d\la\,d\tau \leq CW_k(r,t)^{-1} \cdot \Phi_k(t+r),
\end{equation*}
where $\la_\pm = t-\tau\pm r$, $\Phi_k$ is given by \eqref{Phi} and the constant $C$ is independent of $r, t$.
\end{lemma}

\begin{proof}
Let us first invoke our estimate \eqref{fa} to get
\begin{equation*}
\mathcal{J}_2 \leq \frac{C}{(t-r)^m} \int_0^{t-r} \int_0^{\la_-} \:\frac{(\la+\tau)^m \cdot \la^{1-mp+m}}{\la_+^a
\,(\la_- -\la)^{1-a}} \cdot W_k(\la,\tau)^{-p} \:d\la\,d\tau.
\end{equation*}

\Case{1} When $t\leq 2r$, we use the fact that $\la_+ \geq r$ within the region of integration.  Recalling our definition
\eqref{W1}, we change variables by $x= \la-\tau$ and $y= \la+\tau$ to find
\begin{equation*}
\mathcal{J}_2 \leq \frac{C}{r^a (t-r)^m} \int_0^{t-r} \frac{y^m \br{y}^{-\mu p}}{(t-r-y)^{1-a}} \int_{-y}^y
(x+y)^{1-mp+m} \cdot \br{x}^{-\nu p} \left( 1+ \ln \frac{\br{y}}{\br{x}} \right)^{p\de_{k,m+a}} dx\,dy.
\end{equation*}
Once we now employ Lemma \ref{fa2} to treat the inner integral, we obtain
\begin{equation}\label{In4}
\mathcal{J}_2 \leq \frac{C\Phi_k(t+r)}{r^a (t-r)^m} \int_0^{t-r} \frac{y^m \br{y}^{m-k} \:dy}{(t-r-y)^{1-a}} =
\frac{C\Phi_k(t+r)}{r^a (t-r)^m} \int_0^{t-r} \frac{y^m \br{y}^{-\mu-\nu} \:dy}{(t-r-y)^{1-a}}
\end{equation}
since $\mu+\nu = k-m$.  In view of \eqref{In3} and \eqref{ne}, the desired estimate follows when $t\leq 2r$.

\Case{2} When $t\geq 2r$, it is convenient to express the given integral as the sum of
\begin{equation*}
\mathcal{J}_{21} = \int_0^{2(t-r)/3} \int_0^{\la_-} \frac{\la^{1-mp+2m} \cdot W_k(\la,\tau)^{-p}}{\la_-^m \la_+^a
\,(\la_- -\la)^{1-a}} \:\:d\la\,d\tau
\end{equation*}
and
\begin{equation*}
\mathcal{J}_{22} = \int_{2(t-r)/3}^{t-r} \int_0^{\la_-} \frac{\la^{1-mp+2m} \cdot W_k(\la,\tau)^{-p}}{\la_-^m \la_+^a
\,(\la_- -\la)^{1-a}} \:\:d\la\,d\tau.
\end{equation*}
To treat $\mathcal{J}_{21}$, we proceed as in Case 1 using the inequality $\la_+\geq (t-r)/3$ instead of our previous
$\la_+ \geq r$. Analogously to \eqref{In4}, we now establish
\begin{equation*}
\mathcal{J}_{21} \leq \frac{C\Phi_k(t+r)}{(t-r)^{m+a}} \int_0^{t-r} \frac{y^m \br{y}^{-\mu -\nu}\:dy}{(t-r-y)^{1-a}} \leq
C\Phi_k(t+r) \cdot \br{t-r}^{-\mu-\nu}
\end{equation*}
because of \eqref{In3}.  This does imply the desired estimate whenever $t\geq 2r$.

To treat $\mathcal{J}_{22}$, we use the inequality $\la_+\geq \la_-$ to first obtain
\begin{equation*}
\mathcal{J}_{22} \leq \int_{2(t-r)/3}^{t-r} \int_0^{\la_-} \frac{\la^{1-mp+2m} \cdot W_k(\la,\tau)^{-p}}{\la_-^{m+a}
\,(\la_- -\la)^{1-a}} \:\:d\la\,d\tau.
\end{equation*}
As in Lemma \ref{I2}, $\tau\pm \la$ and $t-r$ are all equivalent here and this makes $W_k(\la,\tau)$ equivalent to
$\br{t-r}^{\mu+\nu}= \br{t-r}^{k-m}$ within the region of integration.  In particular,
\begin{equation*}
\mathcal{J}_{22} \leq C\br{t-r}^{mp-kp} \int_0^{t-r} \int_0^{\la_-} \frac{\la^{1-mp+2m}}{\la_-^{m+a} \,(\la_-
-\la)^{1-a}} \:\:d\la\,d\tau
\end{equation*}
and our next lemma allows us to conclude that
\begin{equation*}
\mathcal{J}_{22} \leq C\br{t-r}^{2-kp+m} = C\br{t-r}^{2-k(p-1)} \cdot \br{t-r}^{m-k}.
\end{equation*}
Recalling the definition of $\Phi_k$, the desired estimate now follows since $m-k=-\mu-\nu$.
\end{proof}

\begin{lemma}\label{J2s}
Let $m\geq a>0$ and $p>1$.  Assuming \eqref{con}, one has
\begin{equation*}
\mathbb{J} \equiv \int_0^{t-r} \int_0^{\la_-} \frac{\la^{1-mp+2m} \:\:d\la\,d\tau}{\la_-^{m+a} \,(\la_- -\la)^{1-a}} \leq
C \br{t-r}^{2-mp+m},
\end{equation*}
where $\la_- = t-r-\tau$ and the constant $C$ depends solely on $a$, $m$ and $p$.
\end{lemma}

\begin{proof}
Since $m$ is positive, we clearly have
\begin{equation*}
\mathbb{J}\leq \int_0^{t-r} \int_0^{\la_-} \frac{\la^{1-mp+m} \:\:d\la\,d\tau}{\la_-^a \,(\la_- -\la)^{1-a}}
\end{equation*}
and the change of variables $y= \la_-= t-r-\tau$ allows us to write
\begin{equation*}
\mathbb{J} \leq \int_0^{t-r}\int_{y/2}^y \frac{\la^{1-mp+m} \:d\la\,dy}{y^a \,(y-\la)^{1-a}} + \int_0^{t-r}\int_0^{y/2}
\frac{\la^{1-mp+m} \:d\la\,dy}{y^a \,(y-\la)^{1-a}} \:.
\end{equation*}
Here, $\la$ is equivalent to $y$ within the former integral and $y-\la$ is equivalent to $y$ within the latter, so we
find that
\begin{equation*}
\mathbb{J} \leq C\int_0^{t-r} y^{1-mp+m-a} \int_{y/2}^y (y-\la)^{a-1} \:d\la\,dy + C\int_0^{t-r} y^{-1} \int_0^{y/2}
\la^{1-mp+m} \:d\la\,dy.
\end{equation*}
Since $a>0$ by assumption and since $2-mp+m>0$ by \eqref{con}, we then obtain
\begin{equation*}
\mathbb{J} \leq C\int_0^{t-r} y^{1-mp+m} \:dy \leq C\br{t-r}^{2-mp+m}.
\end{equation*}
This already establishes the desired estimate, so the proof of our lemma is complete.
\end{proof}

\begin{lemma}\label{Ipm}
Let $W_k$ be as in \eqref{W1} and $\kappa>2$.  Assuming \eqref{ks}, one has
\begin{equation}\label{Ipm1}
\mathcal{I}_\pm \equiv \int_{\max{(t-2r,0)}}^t |\la_\pm|^a \cdot \br{\la_\pm}^{1 -\kappa} \cdot W_k(|\la_\pm|, \tau)^{-1}
\:d\tau \leq Cr^a W_k(r,t)^{-1},
\end{equation}
where $\la_\pm = t-\tau \pm r$ and the constant $C$ is independent of $r,t$.
\end{lemma}

\begin{proof}
Before we turn our attention to the given integrals, let us first check that
\begin{equation}\label{nn}
|\la_\pm|^a \, \br{|\la_\pm| + \tau}^{-\mu} \leq Cr^a \br{t+r}^{-\mu}\quad\quad\text{when $\max(t-2r,0)\leq \tau\leq t.$}
\end{equation}

\Case{1} If either $t\geq 2r$ or $r\leq 1$, then each of $\br{|\la_\pm| + \tau}$ is equivalent to $\br{t+r}$ because
\begin{equation*}
t - r = (t-\tau-r) + \tau \leq |\la_\pm| + \tau \leq t+r
\end{equation*}
are all equivalent when $t\geq 2r$ and all bounded when $t\leq 2r\leq 2$.  Once we now note that
\begin{equation*}
|\la_\pm|^a = |t-\tau \pm r|^a \leq (3r)^a
\end{equation*}
whenever $\tau\geq t-2r$ and $a>0$, our assertion \eqref{nn} follows.

\Case{2} Suppose now that $t\leq 2r$ and $r\geq 1$.  Since $a>0$ and $\mu\leq a$ by \eqref{ks}, we have
\begin{equation*}
|\la_\pm|^a \, \br{|\la_\pm| + \tau}^{-\mu} \leq \br{|\la_\pm| + \tau}^{a-\mu} \leq \br{t+r}^{a-\mu}.
\end{equation*}
Moreover, $r$ and $\br{t+r}$ are equivalent when $r\geq \max(t/2,1)$, so this also implies \eqref{nn}.

Next, we focus on the integrals \eqref{Ipm1}.  Employing our estimate \eqref{nn}, we find
\begin{equation}\label{Ipm2}
\mathcal{I}_\pm \leq Cr^a \br{t+r}^{-\mu} \cdot \mathbb{I}_\pm, \quad\quad \mathbb{I}_\pm \equiv \int_0^t \br{\la_\pm}^{1
-\kappa} \,\br{|\la_\pm|+\tau}^\mu \cdot W_k(|\la_\pm|, \tau)^{-1} \:d\tau.
\end{equation}
Since $\la_\pm = t\pm r-\tau$ by definition, one clearly has
\begin{align*}
\mathbb{I}_\pm
&\leq \int_0^{\max(t\pm r,0)} \br{t\pm r - \tau}^{1 -\kappa} \cdot \br{t\pm r - 2\tau}^{-\nu}
\left( 1+ \ln \frac{\br{t\pm r}}{\br{t\pm r -2\tau}} \right)^{\de_{k,m+a}} d\tau \notag \\
&\quad + \br{t-r}^{-\nu} \int_{\max(t-r,0)}^t \br{\tau-t+r}^{1 -\kappa} \left( 1+ \ln \frac{\br{2\tau-t+r}}{\br{t-r}}
\right)^{\de_{k,m+a}} d\tau.
\end{align*}
Besides, the substitution $x= t\pm r-2\tau$ allows us to write the former integral as
\begin{equation*}
\frac{1}{2} \int_{-(t\pm r)}^{t\pm r} \br{\frac{x+t\pm r}{2}}^{1 -\kappa} \cdot \br{x}^{-\nu} \left( 1+ \ln
\frac{\br{t\pm r}}{\br{x}} \right)^{\de_{k,m+a}} dx.
\end{equation*}
Once we now treat this integral using Lemma \ref{fa1} with $y= t\pm r$, we arrive at
\begin{equation*}
\mathbb{I}_\pm \leq C\br{t\pm r}^{-\nu} + \br{t-r}^{-\nu} \int_{\max(t-r,0)}^t \br{\tau-t+r}^{1 -\kappa} \left( 1+ \ln
\frac{\br{2\tau-t+r}}{\br{t-r}} \right)^{\de_{k,m+a}} d\tau.
\end{equation*}
Here, $0\leq 2\tau-t+r\leq t+r$ within the region of integration and $\kappa>2$, so we find
\begin{equation*}
\mathbb{I}_\pm \leq C\br{t\pm r}^{-\nu} + C\br{t-r}^{-\nu} \left( 1+ \ln \frac{\br{t+r}}{\br{t-r}} \right)^{\de_{k,m+a}}.
\end{equation*}
Since $\nu\geq 0$ by \eqref{ks}, this also implies
\begin{equation*}
\mathbb{I}_\pm \leq C\br{t-r}^{-\nu} \left( 1+ \ln \frac{\br{t+r}}{\br{t-r}} \right)^{\de_{k,m+a}}.
\end{equation*}
Combining the last equation with \eqref{Ipm2}, we thus deduce the desired estimate for $\mathcal{I}_\pm$.
\end{proof}

\begin{lemma}\label{Jpm}
Let $W_k$ be as in \eqref{W1}.  Assume \eqref{k1}, \eqref{con} and \eqref{ks}.  When $p>1$, one then has
\begin{equation*}
\mathcal{J}_\pm \equiv \int_{\max{(t-2r,0)}}^t |\la_\pm|^{a+1-mp+m} \br{\la_\pm}^{-1} \cdot W_k(|\la_\pm|, \tau)^{-p}
\:d\tau \leq Cr^a W_k(r,t)^{-1} \cdot \Phi_k(t+r),
\end{equation*}
where $\la_\pm= t-\tau \pm r$, $\Phi_k$ is given by \eqref{Phi} and the constant $C$ is independent of $r,t$.
\end{lemma}

\begin{proof}
In view of \eqref{nn}, the desired estimate will follow once we show that each of
\begin{equation*}
\mathbb{J}_\pm \equiv \int_0^t |\la_\pm|^{1-mp+m} \br{\la_\pm}^{-1} \br{|\la_\pm| + \tau}^\mu \cdot W_k(|\la_\pm|,
\tau)^{-p} \:d\tau
\end{equation*}
satisfies an inequality of the form
\begin{equation}\label{pr}
\mathbb{J}_\pm \leq C\br{t-r}^{-\nu} \cdot \Phi_k(t+r).
\end{equation}
Let us then proceed as in the previous lemma.  Since $\la_\pm = t\pm r-\tau$, we clearly have
\begin{align}\label{Jpm1}
\mathbb{J}_\pm
&\leq \int_0^{\max(t\pm r,0)} \la_\pm^{1-mp+m} \br{\la_\pm + \tau}^\mu \cdot W_k(\la_\pm, \tau)^{-p} \:d\tau \notag\\
&\quad + \int_{\max(t-r,0)}^t \: (-\la_-)^{1-mp+m} \cdot \br{-\la_-}^{-1} \br{\tau-\la_-}^\mu \cdot W_k(-\la_-,
\tau)^{-p} \:d\tau \notag\\
&\equiv \mathbb{J}_\pm' + \mathbb{J}_-''
\end{align}
and we shall first focus on $\mathbb{J}_\pm'$.  Explicitly, this term is given by
\begin{equation*}
\mathbb{J}_\pm' = \br{t\pm r}^{\mu-\mu p} \int_0^{\max(t\pm r,0)} \la_\pm^{1 -mp+m} \br{\la_\pm -\tau}^{-\nu p} \left( 1+
\ln \frac{\br{t\pm r}}{\br{\la_\pm-\tau}} \right)^{p\de_{k,m+a}} d\tau.
\end{equation*}
Using the substitution $x= \la_\pm -\tau= t\pm r-2\tau$, we may thus express it in the form
\begin{equation*}
\mathbb{J}_\pm' = C\br{t\pm r}^{\mu- \mu p} \int_{-(t\pm r)}^{t\pm r} \: (x+t\pm r)^{1 -mp+m} \br{x}^{-\nu p} \left( 1+
\ln \frac{\br{t\pm r}}{\br{x}} \right)^{p\de_{k,m+a}} dx.
\end{equation*}
Once we now employ Lemma \ref{fa2} with $y= t\pm r$, we find
\begin{equation*}
\mathbb{J}_\pm' \leq C\Phi_k(t\pm r) \cdot \br{t\pm r}^{\mu + m-k} = C\Phi_k(t\pm r) \cdot \br{t\pm r}^{-\nu} \leq
C\Phi_k(t+ r) \cdot \br{t- r}^{-\nu}
\end{equation*}
because $\mu+\nu= k-m$ and $\nu\geq 0$ by \eqref{ks}.  Thus, the desired \eqref{pr} is satisfied by $\mathbb{J}_\pm'$.

To treat the remaining term $\mathbb{J}_-''$, we change variables by $\sigma= -\la_-= \tau-t+r$ and write
\begin{equation}\label{J''}
\mathbb{J}_-'' = \int_{\max(r-t,0)}^r \: \sigma^{1-mp+m} \br{\sigma}^{-1} \br{2\sigma+t-r}^\mu \cdot W_k(\sigma,
\sigma+t-r)^{-p} \:d\sigma.
\end{equation}
Recall the definition \eqref{W1} of our weight function, according to which
\begin{equation}\label{rec}
W_k(\sigma, \sigma+t-r)^{-p} = \br{2\sigma + t-r}^{-\mu p} \br{t-r}^{-\nu p} \left( 1+ \ln
\frac{\br{2\sigma+t-r}}{\br{t-r}} \right)^{p\de_{k,m+a}}
\end{equation}
with $\mu= \min(k-m,a)$ and $\nu= \max(k-m-a,0)$ for some $a>0$.

\Case{1} When $k\leq m$, we have $\mu= k-m\leq 0$ and $\nu=0$, while equation \eqref{J''} reads
\begin{equation*}
\mathbb{J}_-'' = \int_{\max(r-t,0)}^r \: \sigma^{1-mp+m} \br{\sigma}^{-1} \br{2\sigma+t-r}^{-\mu (p-1)} \:d\sigma.
\end{equation*}
Since $0\leq 2\sigma+ t-r\leq t+r$ within the region of integration, this trivially implies
\begin{equation*}
\mathbb{J}_-'' \leq \br{t+r}^{-\mu (p-1)} \int_0^{t+r} \sigma^{1-mp+m} \:d\sigma
\end{equation*}
because $-\mu (p-1)\geq 0$ here.  Using our assumption \eqref{con} that $2-mp+m>0$, we then get
\begin{equation*}
\mathbb{J}_-'' \leq C\br{t+r}^{-\mu (p-1) + 2-mp+m} = C\br{t+r}^{2-k(p-1)} \leq C\Phi_k(t+r)
\end{equation*}
because $\mu= k-m$ here.  This is precisely the desired \eqref{pr}, as $\nu=0$ for this case.

\Case{2} When $k> m$, we have $\mu>0$ and it is convenient to introduce the constant
\begin{equation}\label{de}
\de_* = \left\{ \begin{array}{ccc} \min(\mu(p-1),1/2) &&\text{if\, $k\leq m+a$}\\ 0&&\text{if\, $k> m+a$} \end{array}
\right\}.
\end{equation}
Since $\de_*$ is positive when $k=m+a$, we may then estimate \eqref{rec} as
\begin{align*}
W_k(\sigma, \sigma+t-r)^{-p}
&\leq C\br{2\sigma + t-r}^{\de_* -\mu p} \br{t-r}^{-\nu p} \\
&\leq C\br{2\sigma + t-r}^{\de_* -\mu p} \br{t-r}^{-\nu}
\end{align*}
because $\nu p\geq \nu\geq 0$ by \eqref{ks}.  Inserting this fact in \eqref{J''}, we thus arrive at
\begin{equation*}
\mathbb{J}_-'' \leq C\br{t-r}^{-\nu} \int_{\max(r-t,0)}^r \: \sigma^{1-mp+m} \br{\sigma}^{-1} \br{2\sigma+t-r}^{\de_*
-\mu(p-1)} \:d\sigma.
\end{equation*}
Here, $2\sigma+t-r\geq \sigma$ within the region of integration, so we easily get
\begin{equation*}
\mathbb{J}_-'' \leq C\br{t-r}^{-\nu} \int_0^r \sigma^{1-mp+m} \br{\sigma}^{\de_*-\mu(p-1)-1} \:d\sigma
\end{equation*}
because $\de_\ast\leq \mu(p-1)$ by our choice \eqref{de}.  In particular, we get
\begin{equation*}
\mathbb{J}_-'' \leq C\br{t-r}^{-\nu} \left[ \int_0^{\min(r,1)} \sigma^{1-mp+m} \:d\sigma + \int_{\min(r,1)}^r
\br{\sigma}^{\de_* -(m+\mu)(p-1)} \:d\sigma \right].
\end{equation*}
To deduce the desired estimate \eqref{pr}, it thus suffices to show that
\begin{equation*}
\int_0^{\min(r,1)} \sigma^{1-mp+m} \:d\sigma + \int_{\min(r,1)}^r \br{\sigma}^{\de_* -(m+\mu)(p-1)} \:d\sigma \leq
C\Phi_k(t+r).
\end{equation*}
According to \eqref{con}, the former integral is finite, so it certainly satisfies the last inequality.  Let us then
worry about the latter integral and seek an estimate of the form
\begin{equation}\label{pr2}
\mathbb{K} \equiv \int_0^r \br{\sigma}^{\de_* -(m+\mu)(p-1)} \:d\sigma \leq C\Phi_k(t+r).
\end{equation}

\Subcase{2a} When $k\leq m+a$, we have $\mu= k-m$ and we easily get
\begin{equation*}
\mathbb{K} = \int_0^r \br{\sigma}^{\de_* -k(p-1)} \:d\sigma \leq C\Phi_k(t+r) \int_0^r \br{\sigma}^{\de_*-2} \:d\sigma
\leq C\Phi_k(t+r)
\end{equation*}
because $\de_* < 1$ by our choice \eqref{de}.

\Subcase{2b} When $k> m+a$, we have $\mu=a$ and $\de_*=0$ so that
\begin{equation*}
\mathbb{K} = \int_0^r \br{\sigma}^{-(m+a)(p-1)} \:d\sigma.
\end{equation*}
In view of our assumption \eqref{k1} that $k< (m+a)p-1$, however, we also have
\begin{equation*}
(m+a)(p-1) > k+ 1 -m-a > 1
\end{equation*}
for this subcase, so the last integral is finite.  In particular, the desired \eqref{pr2} follows.
\end{proof}

We are finally in a position to give the

\begin{proof_of}{Theorem \ref{bes}}
To establish our two assertions, we apply Proposition \ref{Les}.  Given a function $G\in \mathcal{C}^1(\Omega_T)$ that
satisfies the singularity condition
\begin{equation}\label{sg3}
G(\la,\tau) = O(\la^{-2m-2+\de}) \quad\quad \text{as\, $\la\to 0$}
\end{equation}
for some fixed $\de>0$, the lemma ensures that
\begin{align}\label{hg}
|D^\b [\mathscr{L}G]|
&\leq Cr^{j-|\b|-m-a} \int_0^t \int_{|\la_-|}^{\la_+} \frac{\la^{m-j+1}}{(\la-\la_-)^{1-a}} \cdot \sum_{s=0}^j
\la^s \,|\d_\la^s G(\la,\tau)| \:d\la\,d\tau \notag\\
&\quad + Cr^{j-|\b|-m} \int_0^{t-r} \int_0^{\la_-} \frac{\la^{2m-j+1}}{\la_-^m \la_+^a \,(\la_- -\la)^{1-a}} \cdot
\sum_{s=0}^j \la^s \,|\d_\la^s G(\la,\tau)| \:d\la\,d\tau \notag\\
&\quad + Cr^{j-|\b|-m-a} \int_{\max(t-2r,0)}^t  |\la_\pm|^{a+m-j+1} \left[ \sum_{s=0}^{j-1} \la^s \,|\d_\la^s G(\la,\tau)|
\right]_{\la= |\la_\pm|} d\tau
\end{align}
whenever $|\b|\leq j\leq 1$.

First, we apply this fact to $G(\la,\tau)= F(u(\la,\tau))$.  By Lemma \ref{bas} with $j=j_0=0$, we have
\begin{equation*}
F(u(\la,\tau)) = O(\la^{-mp}) = O(\la^{-2m-2+\de_1}) \quad\quad \text{as\:\:$\la\to 0$,}
\end{equation*}
where $\de_1= 2-mp+2m$ is positive by \eqref{con}.  In particular, our estimate \eqref{hg} does hold for the special case
$G= F(u)$. Besides, the sums that appear in the right hand side are those of Lemma \ref{bas}, according to which
\begin{equation*}
\sum_{s=0}^{j_0} \la^s \,|\d_\la^s F(u(\la,\tau))| \leq 2Ap ||u||^p \cdot\la^{j-mp} \br{\la}^{j_0-j} \cdot
W_k(\la,\tau)^{-p}, \quad\quad j_0= j, j-1.
\end{equation*}
Once we now insert this fact in \eqref{hg}, we obtain an estimate of the form
\begin{equation*}
|D^\b [\mathscr{L} F(u)]| \leq C||u||^p \cdot r^{j-|\b|-m} \cdot (r^{-a} \mathcal{J}_1 + \mathcal{J}_2 + r^{-a}
\mathcal{J}_\pm),
\end{equation*}
where $\mathcal{J}_1$, $\mathcal{J}_2$ and $\mathcal{J}_\pm$ are as in Lemmas \ref{J1}, \ref{J2} and \ref{Jpm},
respectively.  The assumptions we imposed in these lemmas are not different from the ones imposed in this theorem, except
for the inequality $a\leq m$ that appears in Lemma \ref{J2}.  Nevertheless, our definition \eqref{am} shows that $a\leq
1\leq m$ whenever $n\geq 4$, so we may employ Lemmas \ref{J1}, \ref{J2} and \ref{Jpm} to arrive at
\begin{equation}\label{an1}
|D^\b [\mathscr{L} F(u)]| \leq C||u||^p \cdot r^{j-|\b|-m} \cdot W_k(r,t)^{-1} \cdot \Phi_k(t+r), \quad\quad |\b|\leq
j\leq 1.
\end{equation}
We now claim that this also implies our first assertion \eqref{mb1}, namely that
\begin{equation*}
|D^\b [\mathscr{L} F(u)]| \leq C||u||^p \cdot r^{1-|\b|-m} \br{r}^{|\b|-1} \cdot W_k(r,t)^{-1} \cdot \Phi_k(t+r),
\quad\quad |\b|\leq 1.
\end{equation*}
Indeed, if $r\leq 1$, one may obtain the last inequality through the special case $j=1$ of \eqref{an1}.  If $r\geq 1$, on
the other hand, one may obtain it through the special case $j=|\b|$.

The proof of our assertion \eqref{mb2} regarding the potential term is similar, so we only give a sketch of the proof. In
this case, it suffices to show that
\begin{equation}\label{an2}
|D^\b [\mathscr{L} (Vu)]| \leq CV_0 ||u|| \cdot r^{j-|\b|-m} \cdot W_k(r,t)^{-1}, \quad\quad |\b|\leq j\leq 1.
\end{equation}
We now apply \eqref{hg} with $G(\la,\tau) = V(\la)\cdot u(\la,\tau)$.  Using Lemma \ref{bas} with $j=0$, one easily
checks that the singularity condition \eqref{sg3} holds, and this validates our estimate \eqref{hg} for the special case
$G= Vu$. To treat all three sums that appear in the right hand side, we use the inequality provided by Lemma \ref{bas},
thus arriving at
\begin{equation*}
|D^\b [\mathscr{L} (Vu)]| \leq CV_0 ||u|| \cdot r^{j-|\b|-m} \cdot (r^{-a} \mathcal{I}_1 + \mathcal{I}_2 + r^{-a}
\mathcal{I}_\pm).
\end{equation*}
Here, $\mathcal{I}_1$, $\mathcal{I}_2$ and $\mathcal{I}_\pm$ are given by Lemmas \ref{I1}, \ref{I2} and \ref{Ipm},
respectively. These lemmas are all applicable, as before, so we may invoke them to deduce the desired estimate
\eqref{an2}.
\end{proof_of}

\section{Existence of solutions}\label{bex}
In this section, we give the proof of our existence result, Theorem \ref{et}.  Our first step is to refine the \textit{a
priori} estimates of the previous section, treating the radial derivatives of the Riemann operator in a more efficient
manner.  Since the time derivatives do not appear in our norm \eqref{N}, those are not as important.  We shall merely
need to control them in order to prove uniqueness of solutions using a standard energy argument.

\begin{lemma}[Radial derivatives]\label{ra1}
Fix an integer $n\geq 4$ and let $L$ be the Riemann operator of Lemma \ref{hs}.  Suppose that $f\in \mathcal{C}^1(\R_+)$.
When $r\geq 2t>0$, one then has
\begin{align}\label{Lex2}
|\d_r^{j_0} [Lf](r,t)| \leq Cr^{-j_0} \int_{r-t}^{r+t} |f(\la)| \:d\la  + j_0 C \cdot tr^{-j_0} \sup_{r-t\leq \la\leq
r+t} \: \sum_{s=0}^{j_0} \left[ \la^s \,|f^{(s)}(\la)| \right]
\end{align}
for $j_0= 0,1$ and some constant $C$ that is independent of $r,t$.
\end{lemma}

\begin{proof}
We merely concern ourselves with the case that $n$ is even, as the argument becomes much simpler when $n$ is odd.  Since
$r\geq 2t$, the Riemann operator \eqref{Le} takes the form
\begin{equation}\label{Lf2}
[Lf](r,t) = \frac{1}{2r^{(n-1)/2}} \int_{r-t}^{r+t} \la^{(n-1)/2} f(\la) \cdot U_m(z(\la,r,t)) \:d\la,
\end{equation}
where $U_m$ is given by \eqref{U} and $z(\la,r,t)$ is the rational function \eqref{z}. One may easily check that $0\leq
z(\la,r,t)\leq 1$ whenever $0\leq r-t\leq \la\leq r+t$ and that $z(r\pm t,r,t)= 1$.  Within the region of integration, in
particular, we may express the function of \eqref{U} as
\begin{equation}\label{U2}
U_m(z) = \frac{\sqrt 2}{\pi} \int_z^1 \frac{1}{\sqrt{\sigma - z}} \cdot \frac{T_m(\sigma)}{\sqrt{1- \sigma^2}}
\:\:d\sigma = \frac{\sqrt 2}{\pi} \int_0^1 \frac{1}{\sqrt{\nu (1-\nu)}} \cdot \frac{T_m(\sigma)}{\sqrt{1+ \sigma}}
\:\:d\nu
\end{equation}
by means of the substitution $\sigma = z+ \nu(1-z)$.  Since $z(r\pm t,r,t)= 1$, we then get
\begin{equation}\label{U3}
U_m(z(r\pm t,r,t)) = U_m(1) = \frac{\sqrt 2}{\pi} \int_0^1 \frac{1}{\sqrt{\nu (1-\nu)}} \cdot \frac{T_m(1)}{\sqrt{2}}
\:\:d\nu = T_m(1).
\end{equation}
Using this fact, we now differentiate \eqref{Lf2} to find that
\begin{align*}
\d_r [Lf](r,t)
&= -\frac{n-1}{4r^{(n+1)/2}} \int_{r-t}^{r+t} \la^{(n-1)/2} f(\la) \cdot U_m(z(\la,r,t)) \:d\la \notag\\
&\quad + \frac{1}{2r^{(n-1)/2}} \int_{r-t}^{r+t} \la^{(n-1)/2} f(\la) \cdot \d_r U_m(z(\la,r,t)) \:d\la \notag\\
&\quad + \frac{T_m(1)}{2r^{(n-1)/2}} \cdot \Bigl[ (r+t)^{(n-1)/2} f(r+t) - (r-t)^{(n-1)/2} f(r-t) \Bigr].
\end{align*}
Estimating the last equation and \eqref{Lf2} at the same time, we thus obtain
\begin{align*}
|\d_r^{j_0} [Lf](r,t)| &\leq \sum_{s=0}^{j_0} Cr^{s-j_0-(n-1)/2} \int_{r-t}^{r+t} \la^{(n-1)/2} \,|f(\la)| \cdot
|\d_r^s U_m(z(\la,r,t))| \:d\la \notag \\
&\quad +j_0Cr^{-(n-1)/2} \cdot \Bigl| (r+t)^{(n-1)/2} f(r+t) - (r-t)^{(n-1)/2} f(r-t) \Bigr| \notag\\
&\equiv A_1 + A_2
\end{align*}
for $j_0=0,1$.  When it comes to the boundary terms $A_2$, the mean value theorem yields
\begin{equation*}
A_2 \leq j_0 C\cdot tr^{-(n-1)/2} \sup_{r-t\leq \la\leq r+t} \: \sum_{s=0}^1 \left[ \la^{(n-1)/2-1+s} \,|f^{(s)}(\la)|
\right].
\end{equation*}
Since $r\geq 2t$ by assumption, each $r-t\leq \la\leq r+t$ is equivalent to $r$, so the desired \eqref{Lex2} does hold
for these terms.  When it comes to the integral term $A_1$, we similarly have
\begin{equation*}
A_1\leq Cr^{-j_0} \int_{r-t}^{r+t} |f(\la)| \cdot \sum_{s=0}^{j_0} \la^s \,|\d_r^s U_m(z(\la,r,t))| \:d\la,
\end{equation*}
so the desired \eqref{Lex2} will follow once we know that
\begin{equation}\label{U4}
\sum_{s=0}^1 \la^s \,|\d_r^s U_m(z(\la,r,t))| \leq C, \quad\quad r-t\leq \la\leq r+t.
\end{equation}

Now, a direct differentiation of \eqref{U2} allows us to write
\begin{equation*}
\d_r^s U_m(z(\la,r,t)) = \frac{\sqrt 2}{\pi} \int_0^1 \frac{1}{\sqrt{\nu (1-\nu)}} \cdot \d_r^s \left(
\frac{T_m(\sigma)}{\sqrt{1+ \sigma}} \right) \:d\nu,
\end{equation*}
where $\sigma = z+ \nu(1-z)$ and $z= z(\la,r,t)$.  Since $0\leq z\leq 1$ whenever $0\leq r-t\leq \la\leq r+t$, it is
clear that $0\leq \sigma\leq 1$ within the region of integration.  Thus, one easily finds
\begin{equation}\label{U5}
|\d_r^s U_m(z(\la,r,t))| \leq C|\d_r^s z(\la,r,t)|, \quad\quad r-t\leq \la\leq r+t.
\end{equation}
Here, the rational function $z(\la,r,t)$ is defined by \eqref{z} and satisfies
\begin{equation*}
z(\la,r,t) = \frac{\la^2+r^2- t^2}{2r\la} \:,\quad\quad \d_r z(\la,r,t) = \frac{r^2 - \la^2 + t^2}{2r^2\la} \:.
\end{equation*}
Since $r\geq 2t$ by assumption, we also have $\la,r,t\leq 2r$ and a rather crude estimate gives
\begin{equation*}
|\d_r^s z(\la,r,t)| \leq \frac{Cr^{1-s}}{\la} \:,\quad\quad s=0,1.
\end{equation*}
As $\la$ and $r$ are equivalent by above, we may combine this with \eqref{U5} to deduce \eqref{U4}.
\end{proof}

\begin{corollary}\label{np}
Under the assumptions of Theorem \ref{bes}, one actually has
\begin{equation*}
||\mathscr{L}F(u)|| \leq C_1||u||^p \cdot \Phi_k(T), \quad\quad ||\mathscr{L}(Vu)|| \leq C_1 V_0 ||u||
\end{equation*}
for some constant $C_1$ that is independent of $u$.
\end{corollary}

\begin{proof}
In view of the definition \eqref{N} of our norm, we have to show that
\begin{equation}\label{f1}
|\d_r^{j_0} [\mathscr{L} F(u)]| \leq C||u||^p \cdot r^{1-j_0-m} \br{r}^{j_0-1} \cdot W_k(r,t)^{-1}\cdot \Phi_k(T),
\quad\quad j_0=0,1
\end{equation}
for each $(r,t)\in \Omega_T$, as well as
\begin{equation*}
|\d_r^{j_0} [\mathscr{L} (Vu)]| \leq CV_0 ||u|| \cdot r^{1-j_0-m} \br{r}^{j_0-1} \cdot W_k(r,t)^{-1}, \quad\quad j_0=0,1.
\end{equation*}
Since the latter inequality is a special case of Theorem \ref{bes}, we need only worry about the former.  As another
special case of Theorem \ref{bes}, we do have the estimate
\begin{equation*}
|\d_r^{j_0} [\mathscr{L} F(u)]| \leq C||u||^p \cdot r^{1-j_0-m} \br{r}^{j_0-1} \cdot W_k(r,t)^{-1} \cdot \Phi_k(t+r),
\end{equation*}
which implies \eqref{f1} when $\Phi_k(t+r)\leq C\Phi_k(t)$, hence when either $r\leq 2t$ or $r\leq 2$.  In what follows,
we may thus focus on the case $r\geq \max(2t,2)$.  When it comes to
\begin{equation*}
\d_r^{j_0} [\mathscr{L} F(u)](r,t) = \int_0^t \d_r^{j_0} [LF(u(\cdot\,,\tau))](r,t-\tau) \:d\tau,
\end{equation*}
we have $r\geq 2(t-\tau)$ within the region of integration, so Lemma \ref{ra1} applies to give
\begin{align*}
|\d_r^{j_0} [\mathscr{L} F(u)]|
&\leq Cr^{-j_0} \int_0^t \int_{r-t+\tau}^{r+t-\tau} |F(u(\la,\tau))| \:d\la\,d\tau \\
&\quad + Cr^{-j_0} \int_0^t (t-\tau) \sup_{r-t+\tau\leq \la\leq r+t-\tau} \: \sum_{s=0}^{j_0} \Bigl[ \la^s \,|\d_\la^s
F(u(\la,\tau))| \Bigr] \:d\tau.
\end{align*}
According to Lemma \ref{bas} with $j=j_0$, we also have
\begin{equation*}
\sum_{s=0}^{j_0} \la^s |\d_\la^s F(u(\la,\tau))| \leq 2Ap ||u||^p \cdot \la^{j_0-mp} \cdot W_k(\la,\tau)^{-p}, \quad\quad
j_0=0,1
\end{equation*}
so we may combine the last two equations to arrive at
\begin{align*}
|\d_r^{j_0} [\mathscr{L} F(u)]|
&\leq C||u||^p \cdot r^{-j_0} \int_0^t \int_{r-t+\tau}^{r+t-\tau} \la^{j_0-mp} \cdot W_k(\la,\tau)^{-p} \:d\la\,d\tau \\
&\quad + C||u||^p \cdot r^{-j_0} \int_0^t (t-\tau) \sup_{r-t+\tau\leq \la\leq r+t-\tau} \la^{j_0-mp} \cdot
W_k(\la,\tau)^{-p}\:d\tau.
\end{align*}
Due to our assumption that $r\geq 2t$, one has $\la\geq r-t+\tau \geq t+\tau \geq 2\tau$ within the region of
integration, whence $W_k(\la,\tau)$ is equivalent to $\br{\la}^{k-m}$.  In addition,
\begin{equation*}
\frac{r}{2} \leq r-t+ \tau \leq \la \leq r+t-\tau \leq \frac{3r}{2}
\end{equation*}
so $\br{\la}$ is equivalent to $r$ itself and we get
\begin{align*}
|\d_r^{j_0} [\mathscr{L} F(u)]|
&\leq C||u||^p \cdot r^{-j_0} \int_0^t \int_{r-t+\tau}^{r+t-\tau} r^{j_0-mp} \cdot r^{mp-kp} \:d\la\,d\tau \\
&\quad + C||u||^p \cdot r^{-j_0} \int_0^t (t-\tau) \cdot r^{j_0-mp} \cdot r^{mp-kp} \:d\tau.
\end{align*}
Since $k(p-1)\geq 0$ by assumption, the last equation easily leads to
\begin{align*}
|\d_r^{j_0} [\mathscr{L} F(u)]| \leq C||u||^p \cdot r^{-kp} \br{t}^2 \leq C||u||^p \cdot r^{-k} \cdot \br{t}^{2-k(p-1)}
\end{align*}
whenever $r\geq \max(2t,2)$.  Besides, $r$, $\br{r}$ and $\br{r\pm t}$ are all equivalent here, hence
\begin{align*}
|\d_r^{j_0} [\mathscr{L} F(u)]| \leq C||u||^p \cdot r^{-k} \cdot \Phi_k(t) \leq C||u||^p \cdot r^{-m} W_k(r,t)^{-1} \cdot
\Phi_k(t)
\end{align*}
by the definition of $\Phi_k$.  This does imply the desired estimate \eqref{f1} whenever $r\geq 2$.
\end{proof}

Following Kubo \cite{Kb1}, we shall now introduce the auxiliary norm
\begin{equation}\label{Na}
|||u||| = \sup_{(r,t)\in \Omega_T} |u(r,t)|\cdot r^m W_k(r,t).
\end{equation}
Since $r\leq \br{r}$, a comparison with our previous norm \eqref{N} gives $|||u|||\leq ||u||$.

\begin{lemma}\label{bas2}
Let $p>1$ and $X$ be the Banach space \eqref{X}.  With $u,v\in X$ arbitrary, set
\begin{equation}\label{N1}
M(u,v) = |||u - v||| \cdot \Bigl( ||u||^{p-1} + ||v||^{p-1} \Bigr)
\end{equation}
and
\begin{equation}\label{N2}
N(u,v) = ||u-v|| \cdot \Bigl( ||u||^{p-1} + ||v||^{p-1} \Bigr) + |||u-v|||^{p-1} \cdot \Bigl( ||u|| + ||v|| \Bigr).
\end{equation}
Assuming that \eqref{F} holds, one then has
\begin{equation}\label{ma3}
|F(u(\la,\tau))-F(v(\la,\tau))| \leq CM(u,v) \cdot \la^{-mp} \cdot W_k(\la,\tau)^{-p}
\end{equation}
and also
\begin{equation}\label{ma4}
\sum_{s=0}^{j_0} \la^s \,|\d_\la^s [F(u(\la,\tau))-F(v(\la,\tau))]| \leq CN(u,v) \cdot \la^{j-mp} \br{\la}^{j_0-j} \cdot
W_k(\la,\tau)^{-p}
\end{equation}
whenever $0\leq j,j_0\leq 1$ and $(\la,\tau)\in \Omega_T$.  Besides, the constant $C$ is independent of $\la,\tau$.
\end{lemma}

\begin{proof}
Let us first focus on the derivation of \eqref{ma3}. Since
\begin{equation*}
F(u)-F(v) = (u-v) \int_0^1 F'(\theta u + (1-\theta)v) \:d\theta,
\end{equation*}
our assumption \eqref{F} on $F$ easily leads to
\begin{equation}\label{es1}
|F(u)-F(v)| \leq C|u-v| \cdot \Bigl( |u|^{p-1} + |v|^{p-1} \Bigr).
\end{equation}
Using the norm \eqref{Na} for $u-v$ and the norm \eqref{N} for the other factor, we then get
\begin{equation*}
|F(u)-F(v)| \leq C \,|||u-v||| \cdot \la^{p-1-mp} \br{\la}^{1-p} \cdot W_k(\la,\tau)^{-p} \cdot \Bigl( ||u||^{p-1} +
||v||^{p-1} \Bigr).
\end{equation*}
In view of the definition \eqref{N1} of $M(u,v)$, this does imply \eqref{ma3} whenever $p>1$.

Next, we turn to \eqref{ma4}.  To treat the summand for the index $s=0$, we have to show that
\begin{equation}\label{ma41}
|F(u(\la,\tau))-F(v(\la,\tau))| \leq CN(u,v) \cdot \la^{j-mp} \br{\la}^{j_0-j} \cdot W_k(\la,\tau)^{-p}
\end{equation}
whenever $0\leq j,j_0\leq 1$ and $(\la,\tau)\in \Omega_T$.  Using the norm \eqref{N} to now estimate both factors in the
right hand side of \eqref{es1}, we find
\begin{equation*}
|F(u)-F(v)| \leq C \,||u-v|| \cdot \la^{p-mp} \br{\la}^{-p} \cdot W_k(\la,\tau)^{-p} \cdot \Bigl( ||u||^{p-1} +
||v||^{p-1} \Bigr).
\end{equation*}
Moreover, we have $j\leq 1< p$ by assumption, so the definition \eqref{N2} of $N(u,v)$ gives
\begin{equation*}
|F(u)-F(v)| \leq CN(u,v) \cdot \la^{j-mp} \br{\la}^{-j} \cdot W_k(\la,\tau)^{-p}.
\end{equation*}
This also implies the desired estimate \eqref{ma41} because $j_0\geq 0$.

To finish the proof of \eqref{ma4}, it remains to treat the summand for the index $s=1$.  Since this summand is only
present when $j_0=1$, it suffices to show that
\begin{equation}\label{ma42}
\la \,|\d_\la [F(u(\la,\tau))-F(v(\la,\tau))]| \leq CN(u,v) \cdot \la^{j-mp} \br{\la}^{1-j} \cdot W_k(\la,\tau)^{-p}
\end{equation}
whenever $j=0,1$ and $(\la,\tau)\in \Omega_T$.  Now, one clearly has
\begin{equation*}
|\d_\la [F(u)-F(v)]| \leq |F'(u)-F'(v)|\cdot |\d_\la u| + |\d_\la (u-v)| \cdot |F'(v)|
\end{equation*}
by the triangle inequality, so our assumption \eqref{F} on $F$ leads to
\begin{equation*}
|\d_\la [F(u)-F(v)]| \leq Ap \,|u-v|^{p-1} \cdot |\d_\la u| + |\d_\la (u-v)| \cdot Ap \,|v|^{p-1}.
\end{equation*}
Using the norm \eqref{Na} for $u-v$ and the norm \eqref{N} for all the other factors, we then get
\begin{align*}
|\d_\la [F(u)-F(v)]|
&\leq Ap\, |||u-v|||^{p-1} \cdot ||u|| \cdot \la^{-mp} \cdot W_k(\la,\tau)^{-p} \\
&\quad + Ap\, ||u-v|| \cdot ||v||^{p-1} \cdot \la^{p-1-mp} \br{\la}^{1-p} \cdot W_k(\la,\tau)^{-p}.
\end{align*}
In view of the definition \eqref{N2} of $N(u,v)$, this actually implies
\begin{equation*}
\la \,|\d_\la [F(u)-F(v)]| \leq Ap \cdot N(u,v) \cdot \Bigl( \la^{1-mp} + \la^{p-mp} \br{\la}^{1-p} \Bigr) \cdot
W_k(\la,\tau)^{-p}.
\end{equation*}
Since we also have $j\leq 1 < p$ by assumption, we may thus conclude that
\begin{equation*}
\la \,|\d_\la [F(u)-F(v)]| \leq 2Ap \cdot N(u,v) \cdot \la^{j-mp} \br{\la}^{1-j} \cdot W_k(\la,\tau)^{-p}.
\end{equation*}
This is precisely the desired estimate \eqref{ma42}, so the proof of our lemma is complete.
\end{proof}

\begin{corollary}\label{la}
Let $u,v\in X$.  Define $M(u,v)$ and $N(u,v)$ by \eqref{N1} and \eqref{N2}, respectively.  Under the assumptions of
Theorem \ref{bes}, one then has
\begin{subequations}\label{la1}
\begin{align}
||\mathscr{L}(F(u) - F(v))||
&\leq C_2 N(u,v) \cdot \Phi_k(T), \label{la11}\\
||\mathscr{L}(Vu -Vv)|| &\leq C_2V_0 ||u-v|| \label{la12}
\end{align}
\end{subequations}
as well as
\begin{subequations}\label{la2}
\begin{align}
|||\mathscr{L}(F(u) - F(v))||| &\leq C_2 M(u,v) \cdot \Phi_k(T), \label{la21}\\
|||\mathscr{L}(Vu - Vv)||| &\leq C_2V_0 |||u-v|||.\label{la22}
\end{align}
\end{subequations}
Besides, the constant $C_2$ is independent of $u, v$.
\end{corollary}

\begin{proof}
Except for constant factors, our first two assertions \eqref{la1} are the exact analogues of Corollary \ref{np}.  In
fact, \eqref{la12} does follow from Corollary \ref{np}, according to which
\begin{equation*}
||\mathscr{L} (Vu-Vv)|| = ||\mathscr{L} (V(u-v))|| \leq C_1V_0 ||u-v||.
\end{equation*}
As for \eqref{la11}, our previous approach applies almost verbatim.  More precisely, the estimate that Lemma \ref{bas}
provided before was
\begin{equation*}
\sum_{s=0}^{j_0} \la^s \:|\d_\la^s F(u(\la,\tau))| \leq 2Ap ||u||^p \cdot\la^{j-mp} \br{\la}^{j_0-j} \cdot
W_k(\la,\tau)^{-p}
\end{equation*}
for $0\leq j,j_0\leq 1$ and each $(\la,\tau)\in \Omega_T$.  In this case, an analogous estimate \eqref{ma4} is provided
by Lemma \ref{bas2}, so one may establish \eqref{la11} exactly as before.

Our last two assertions \eqref{la2} follow in a similar fashion as well, so we only give a sketch of their proof. Here,
we apply Proposition \ref{Les} with $j=|\b|=0$.  Given a function $G\in \mathcal{C}^1(\Omega_T)$ that satisfies the
singularity condition \eqref{sg2} for some fixed $\de>0$, the lemma ensures that
\begin{align}\label{hg2}
|[\mathscr{L}G](r,t)|
&\leq Cr^{-m-a} \int_0^t \int_{|\la_-|}^{\la_+} \frac{\la^{m+1}}{(\la-\la_-)^{1-a}} \cdot |G(\la,\tau)| \:d\la\,d\tau
\notag\\
&\quad + Cr^{-m} \int_0^{t-r} \int_0^{\la_-} \frac{\la^{2m+1}}{\la_-^m \la_+^a \,(\la_- -\la)^{1-a}} \cdot
|G(\la,\tau)| \:d\la\,d\tau,
\end{align}
where $\la_\pm = t-\tau \pm r$.

First, we take $G(\la,\tau)= F(u(\la,\tau))- F(v(\la,\tau))$ and use \eqref{ma3} to see that the singularity condition
\eqref{sg2} holds in this case.  Once we now estimate $G(\la,\tau)$ using \eqref{ma3}, we get
\begin{align*}
|\mathscr{L}(F(u)-F(v))|
&\leq CM(u,v) \cdot r^{-m-a} \int_0^t \int_{|\la_-|}^{\la_+} \frac{\la^{m+1-mp} \cdot W_k(\la,\tau)^{-p}}{(\la
-\la_-)^{1-a}}  \:\:d\la\,d\tau \\
&\quad + CM(u,v) \cdot r^{-m} \int_0^{t-r} \int_0^{\la_-} \frac{\la^{2m+1-mp} \cdot W_k(\la,\tau)^{-p}}{\la_-^m \la_+^a
\,(\la_- -\la)^{1-a}} \:\:d\la\,d\tau \\
&= CM(u,v) \cdot r^{-m} \cdot (r^{-a} \mathcal{J}_1 + \mathcal{J}_2),
\end{align*}
where $\mathcal{J}_1$ and $\mathcal{J}_2$ are as in Lemmas \ref{J1} and \ref{J2}, respectively.  Thus, we find
\begin{equation*}
|\mathscr{L}(F(u)-F(v))| \leq CM(u,v) \cdot r^{-m} \cdot W_k(r,t)^{-1} \cdot \Phi_k(t+r).
\end{equation*}
In view of the definition \eqref{Na} of our auxiliary norm, this already implies the desired \eqref{la21} when
$\Phi_k(t+r)\leq C\Phi_k(t)$, hence when either $r\leq 2t$ or $r\leq 2$.  When $r\geq \max(2t,2)$, on the other hand, we
may apply Lemma \ref{ra1} with $j_0=0$ to find that
\begin{equation*}
|\mathscr{L}(F(u)-F(v))| \leq C \int_0^t \int_{r-t+\tau}^{r+t-\tau} |F(u(\la,\tau)) - F(v(\la,\tau))| \:d\la\,d\tau.
\end{equation*}
As in the proof of Corollary \ref{np}, $\la$, $\la\pm \tau$, $r\pm t$ and $r$ are all equivalent within the region of
integration here, so one may easily employ \eqref{ma3} to establish \eqref{la21}.

Since the derivation of \eqref{la22} is more straightforward, we are going to omit it.
\end{proof}

We are finally in a position to give the

\begin{proof_of}{Theorem \ref{et}}
Our iteration argument is quite similar to that of \cite{Kb1}, so we only give a sketch of the proof. As we have already
noted, one may decrease the decay rates $k,\kappa$ to ensure that \eqref{k1} through \eqref{ka} hold without loss of
generality. Let $C_0$, $C_1$ and $C_2$ be the constants that appear in Lemma \ref{dec}, Corollary \ref{np} and Corollary
\ref{la}, respectively.  In order to proceed, we assume that $V_0,\e$ are so small that
\begin{equation}\label{ch1}
2C_0 \e < 1,
\end{equation}
\begin{equation}\label{ch2}
4C_i V_0 + 4C_i \cdot (2C_0\e)^{p-1} \cdot \Phi_k(1) < 1, \quad\quad i=1,2.
\end{equation}
Recall that $\Phi_k(t) = \br{t}^{\max(2-k(p-1),0)}$. In the supercritical case $k\geq 2/(p-1)$, this function is
identically equal to $1$, so one has
\begin{equation}\label{ch3}
4C_3 V_0 + 4C_3 \cdot (2C_0\e)^{p-1}\cdot \Phi_k(T) \leq 1, \quad\quad C_3 = \max(C_1,C_2)
\end{equation}
for any $T>0$ by above.  In the subcritical case, on the other hand, the last inequality does hold with equality for some
$T>1$.  For this case, in particular, we take $T>1$ such that
\begin{equation*}
4C_3 V_0 + 4C_3 \cdot (2C_0\e)^{p-1}\cdot \br{T}^{2-k(p-1)} = 1.
\end{equation*}
Due to the equivalence of $\br{T}$ with $T$, we thus obtain the lower bound
\begin{equation*}
T\geq C\e^{-(p-1)/(2-k(p-1))}
\end{equation*}
that \eqref{lbf} asserts for the subcritical case $k< 2/(p-1)$.

The Banach space $X$ of interest was introduced in \eqref{X} and we shall henceforth focus on its subset $X_\de$
consisting of all $u\in X$ with $||u|| < \de$, where
\begin{equation}\label{ch4}
\de = \min \left( 1 \:,\:\left( \frac{1-4C_3V_0}{4C_3\Phi_k(T)} \right)^{1/(p-1)} \right).
\end{equation}
For this particular choice of $\de$, Lemma \ref{dec} easily leads to the estimate $||u_0||\leq C_0\e \leq \de/2$ because
of \eqref{ch1} and \eqref{ch3}.  This means that $u_0\in X_\de$.  Let us then recursively define
\begin{equation}\label{it}
u_{i+1} = u_0 + \mathscr{L} F(u_i) - \mathscr{L}(Vu_i)
\end{equation}
for each $i\geq 0$.  Using Corollary \ref{np} and our choice \eqref{ch4} of $\de$, one easily finds that $u_{i+1}\in
X_\de$ whenever $u_i\in X_\de$.  In particular, the whole sequence $\{u_i\}$ lies in $X_\de$ by induction.

Next, we employ Corollary \ref{la}.  Using its second conclusion \eqref{la2}, we are able to establish the contraction
estimate
\begin{equation}\label{cot}
|||\mathscr{L}(F(u)-F(v))||| + |||\mathscr{L}(Vu-Vv)||| \leq \frac{1}{2} \cdot |||u-v|||, \quad\quad u,v\in X_\de.
\end{equation}
Using repeated applications of this fact and the first conclusion \eqref{la1} of Corollary \ref{la}, we deduce that
$\{u_i\}$ is a Cauchy sequence in $X_\de$, hence also convergent.  Let $u\in X_\de$ be the limit of this sequence. In
view of \eqref{it} then, $u$ satisfies the integral equation
\begin{equation}\label{ite}
u= u_0 + \mathscr{L}F(u) - \mathscr{L}(Vu).
\end{equation}
In view of Lemmas \ref{new} and \ref{hs}, it also satisfies the nonlinear wave equation \eqref{nlf}.

Finally, we prove the uniqueness assertion of our theorem.  According to Lemma \ref{dec} and Theorem \ref{bes}, our
solution \eqref{ite} to \eqref{nlf} is such that $u\in \mathcal{C}^1(\Omega_T)\cap X_\de$ and
\begin{equation}\label{sg}
|\d_r u| + |\d_t u| = O(r^{-m}) = O(r^{-(n-1)/2 + a}) \quad \text{as $r \rightarrow 0$,}
\end{equation}
where $a>0$ is defined by \eqref{am}. Given some other solution $v$ with the same properties,
\begin{equation*}
w =v- \mathscr{L}F(v) + \mathscr{L}(Vv)
\end{equation*}
satisfies the homogeneous wave equation \eqref{he}.  By Lemma \ref{dec} and Theorem \ref{bes}, $w$ is also subject to
\eqref{sg}, so the uniqueness assertion of Lemma 3.2 in \cite{KKo} implies that $w= u_0$. In other words, $v\in X_\de$
satisfies the integral equation \eqref{ite}. Since that equation has at most one solution in $X_\de$ by our contraction
estimate \eqref{cot}, we may conclude that $v=u$.
\end{proof_of}

\section{Blow-up due to the initial data}\label{BD}
In this section, we establish our blow-up result for initial data of subcritical decay rates.  To merely focus on the
behavior of the initial data at infinity, we shall fix a constant $R>0$ and introduce assumptions for the initial data
only when $|x|> R$.  Before we can deal with such data, however, we shall need to refine a result of Glassey \cite{Gl1}
regarding the positivity of the Riemann operator in any space dimension $n\geq 1$; see also \cite{Ra, Ta2}.

\begin{lemma}\label{pos}
Let $L,\mathscr{L}$ denote the Riemann and Duhamel operators of Lemmas \ref{hs} and \ref{new}, respectively. Given any
$R> 0$, the following properties then hold for some constant $\b_n > 1$ that only depends on $n$.

\begin{itemize}
\item[(a)]
Assuming $f\colon (R,\infty)\to \R$ is non-negative and continuous, one has
\begin{equation*}
[Lf](r,t) \geq 0\:\: \text{whenever\, $r> \max(\b_n t, t+R)$ \,and\, $t\geq 0$.}
\end{equation*}

\item[(b)]
Assuming $g\colon \R^2\to \R$ is continuous, let $r> \max(\b_nt, t+R)$ and $t \geq 0$ be now fixed.  In the case that
$g(\la,\tau)\geq 0$ whenever $\la > \max(\b_n\tau, \tau+R)$ and $0\leq \tau \leq t$, one has
\begin{equation*}
[\mathscr{L}g](r,t) \geq 0.
\end{equation*}
\end{itemize}
\end{lemma}

\begin{proof}
We merely concern ourselves with the case that $n$ is even, as the argument becomes simpler when $n$ is odd.  First, we
focus on part (a) and take $r> t+R$ to be arbitrary.  The Riemann operator \eqref{Le} is then of the form
\begin{equation}\label{Lf}
[Lf](r,t) = \frac{1}{2r^{(n-1)/2}} \int_{r-t}^{r+t} \la^{(n-1)/2} f(\la) \cdot U_m(z(\la,r,t))\:d\la,
\end{equation}
where $U_m$ is given by \eqref{U} and $z(\la,r,t)$ is the rational function \eqref{z}.  Since $\la\geq r-t> R$ within the
region of integration, we have $f(\la)\geq 0$ by assumption.  To establish part (a), it thus suffices to find a constant
$\b_m > 1$ such that
\begin{equation}\label{gn}
U_m(z(\la,r,t)) > 0\quad \text{whenever\, $r> \b_m t$ \,and\, $r+t\geq \la\geq r-t$.}
\end{equation}

Now, our computation \eqref{U3} shows that $U_m(1)= T_m(1)$, where $T_m$ is the $m$th Tchebyshev polynomial. Using the
fact that $T_m(1)= 1$, we may then choose some constant $0 < \al_m < 1$ such that $U_m$ is positive on $[\al_m,1]$.
Setting
\begin{equation}\label{gam}
\b_m = \frac{1}{1-\al_m} > 1,
\end{equation}
it thus suffices to show
\begin{equation}\label{gn2}
z(\la,r,t) \in [\al_m, 1]\quad \text{whenever\, $r > \b_m t$ \,and\, $r+t\geq \la\geq r-t$.}
\end{equation}
It is easy to check that the rational function $z(\la,r,t)$ of \eqref{z} satisfies
\begin{equation*}
z(\la,r,t) -1 = \frac{\la^2+r^2-t^2}{2r\la} -1 = \frac{(\la-r+t)(\la-r-t)}{2r\la} \leq 0
\end{equation*}
whenever $r+t\geq \la \geq r-t$.  Assuming that $r> \b_m t$ as well, we now get
\begin{equation*}
z(\la,r,t) = \frac{\la^2+r^2-t^2}{2r\la} \geq \frac{(\la+r+t)(r-t)}{2r\la} \geq \frac{r-t}{r} > \al_m
\end{equation*}
because $r> \b_m t$ if and only if $r-t > \al_m r$.  This proves our claim \eqref{gn2}, so part (a) follows.

Next, we focus on part (b).  Suppose $r> \max(\b_mt, t+R)$ and $t\geq 0$.  In view of \eqref{Lf}, we may then express the
Duhamel operator \eqref{Du} in the form
\begin{equation*}
[\mathscr{L}g](r,t) = \frac{1}{2r^{(n-1)/2}} \int_0^t \int_{r-t+\tau}^{r+t-\tau} \la^{(n-1)/2} g(\la,\tau) \cdot
U_m(z(\la,r,t-\tau)) \:d\la\,d\tau.
\end{equation*}
Since $r> \b_mt$, the rightmost factor in the integrand is positive by \eqref{gn}.  To finish the proof using our
positivity assumption on $g$, it remains to note that
\begin{equation*}
\la\geq r-t+\tau > \max((\b_m-1)t, R) + \tau \geq \max(\b_m \tau, R+\tau)
\end{equation*}
within the region of integration, as $r> \max(\b_m t, t+R)$ by above and $\b_m > 1$ by \eqref{gam}.
\end{proof}

\begin{lemma}[Uniqueness]
Let $R,T>0$ be arbitrary and define the region $\Theta_T$ by
\begin{equation}\label{Theta}
\Theta_T = \{ (r,t)\in \R_+^2 \:\::\:\: r> t+R,\quad 0\leq t< T\}.
\end{equation}
Suppose $g\colon (R,\infty) \times \R\to \R$ is a $\mathcal{C}^1$ function with $g(r,0)\equiv 0$ on $(R,\infty)$ and
$\psi\colon (R,\infty) \to \R$ is continuous. Then, the equation
\begin{equation}\label{nlb}
\left\{
\begin{array}{rll}
\d_t^2 u -\d_r^2 u - \dfrac{n-1}{r} \cdot \d_r u = g(r,u) \quad\quad &\text{in\, $\Theta_T$}\\
u(r,0)= 0; \quad \d_t u(r,0)= \psi(r) \quad\quad &\text{in\, $(R,\infty)$}
\end{array}\right.
\end{equation}
has at most one solution $u$ which is continuous, locally bounded and of locally finite energy in $\Theta_T$, namely,
such that $u\in L^\infty_{\text{loc}}(\Theta_T)$ and $r^{(n-1)/2} (|u|+|\d_r u| + |\d_t u|) \in
L^2_{\text{loc}}(\Theta_T)$.
\end{lemma}

\begin{proof}
If $u_1, u_2$ are two solutions having the desired properties, then $w= u_1- u_2$ satisfies
\begin{equation*}
\d_t^2 w -\d_r^2 w - \dfrac{n-1}{r} \cdot \d_r w = g(r,u_1) - g(r,u_2) \quad\quad \text{in\, $\Theta_T$}
\end{equation*}
and vanishes on $(R,\infty) \times \{ t=0 \}$.  Using a quite standard computation, one also finds
\begin{align*}
\frac{d}{dt} \:\int_{R+t}^{R_0-t} (w_r^2 + w_t^2) \cdot r^{n-1}\,dr
&= - \Bigl[ (w_r - w_t)^2 \cdot r^{n-1} \Bigr]_{r=R_0-t}  - \Bigl[ (w_r + w_t)^2 \cdot r^{n-1} \Bigr]_{r=R+t} \\
&\quad + 2\int_{R+t}^{R_0-t} w_t \cdot [g(r,u_1) - g(r,u_2)] \cdot r^{n-1}\:dr,
\end{align*}
where $w_r= \d_r w$, $w_t= \d_t w$ and $R_0>0$ is arbitrary.  This gives the local energy inequality
\begin{equation*}
\frac{d}{dt} \:\int_{R+t}^{R_0-t} (w_r^2 + w_t^2) \cdot r^{n-1}\,dr \leq C\int_{R+t}^{R_0-t} |w_t w| \cdot r^{n-1}\,dr
\end{equation*}
because $g\in \mathcal{C}^1$ and since $u_1,u_2$ are locally bounded in $\Theta_T$.  Using H\"older's and Gronwall's
inequalities, we deduce that $w= 0$ in $\Theta_T$.  In particular, we deduce that $u_1= u_2$ in $\Theta_T$.
\end{proof}

\begin{lemma}[Existence]\label{ex}
Let $R>0$ be arbitrary and let $g,\psi$ be as in the previous lemma.  Then, there exists some $T>0$ and a unique solution
$u$ to \eqref{nlb} which is continuous, locally bounded and of locally finite energy in $\Theta_T$.
\end{lemma}

\begin{proof}
We merely concern ourselves with the case that $n$ is even, as the argument becomes simpler when $n$ is odd.  Since $r >
t$ within $\Theta_T$, the Riemann operator \eqref{Le} takes the form
\begin{equation}\label{Lf3}
[L\psi](r,t) = \frac{1}{2r^{(n-1)/2}} \int_{r-t}^{r+t} \la^{(n-1)/2} \psi(\la) \cdot U_m(z(\la,r,t)) \:d\la,
\end{equation}
where $U_m$ is given by \eqref{U} and $z(\la,r,t)$ is the rational function \eqref{z}. One may easily check that $0\leq
z(\la,r,t)\leq 1$ whenever $0\leq r-t\leq \la\leq r+t$, so the function \eqref{U} is such that
\begin{equation*}
|U_m(z)| = \left| \frac{\sqrt 2}{\pi} \int_z^1 \frac{1}{\sqrt{\sigma - z}} \cdot \frac{T_m(\sigma)}{\sqrt{1- \sigma^2}}
\:\:d\sigma \right| \leq C(m) \int_z^1 \frac{1}{\sqrt{\sigma-z}} \cdot \frac{1}{\sqrt{1-\sigma}} \:\:d\sigma = C'(m).
\end{equation*}
Estimating \eqref{Lf3} rather crudely, we deduce that
\begin{equation}\label{Lf4}
|[L\psi](r,t)| \leq C(m) \int_{r-t}^{r+t} |\psi(\la)| \:d\la.
\end{equation}
In view of \eqref{Lf3} and \eqref{Lf4} then, $u_0= L\psi$ is both continuous and locally bounded in $\Theta_T$.

With $\mathscr{L}$ the Duhamel operator of Lemma \ref{new}, we now recursively define a sequence $\{u_i\}$ by setting
$u_{i+1} = u_0 + \mathscr{L}g(r,u_i)$ for each $i\geq 0$.  Using \eqref{Lf3} and \eqref{Lf4}, one easily checks that each
$u_i$ is continuous and locally bounded in $\Theta_T$, while
\begin{equation*}
|u_{i+1}(r,t) - u_i(r,t)| \leq C(m) \int_0^t \int_{r-t+\tau}^{r+t-\tau} |g(\la,u_i) - g(\la,u_{i-1})| \:\:d\la\, d\tau.
\end{equation*}
Given a fixed but arbitrary $R_0>0$, we thus arrive at
\begin{equation*}
\sup_{r+t\leq R_0} |u_{i+1} - u_i| \leq C(m, R_0) \cdot T^2 \cdot \sup_{r+t\leq R_0} |u_i - u_{i-1}|.
\end{equation*}
Iterating the last inequality and choosing $T$ to be sufficiently small, we find that $u_i$ converges uniformly to a
function $u$ which is continuous and locally bounded in $\Theta_T$.

As in the proof of Theorem \ref{et}, one easily checks that $u$ also satisfies \eqref{nlb}.  Our assertion that $u$ is of
locally finite energy follows similarly from \eqref{Lf3}, so we shall omit the details.  As for our uniqueness assertion,
this has already been established in our previous lemma.
\end{proof}

\begin{lemma}[Comparison Lemma]\label{cl}
Using the notation and assumptions of the previous two lemmas, assume also that $\d_u g(r,u)\geq 0$ on $(R,\infty) \times
[0,\infty)$.  When $\psi, \widetilde{\psi}\colon (R,\infty)\to \R$ are such that $\psi\geq \widetilde{\psi}\geq 0$, the
corresponding solutions $u, \widetilde{u}$ to \eqref{nlb} are then such that
\begin{equation}\label{Th2}
u\geq \widetilde{u}\geq 0 \quad\text{in\, $\Theta_T' = \{ (r,t)\in \R_+^2 \:\::\:\: r> \max(\b_n t, t+R),\quad 0\leq t<
T\}$.}
\end{equation}
Here, $\b_n$ denotes the constant of Lemma \ref{pos}, while $T$ denotes the lifespan of $u$.
\end{lemma}

\begin{proof}
According to part (a) of Lemma \ref{pos}, $u_0= L\psi$ and $\widetilde{u}_0= L\widetilde{\psi}$ are related by
\begin{equation*}
u_0 \geq \widetilde{u}_0 \geq 0 \quad\text{in\, $\Theta_T'$.}
\end{equation*}
Suppose that we know $u_i \geq \widetilde{u}_i \geq 0$ in $\Theta_T'$ for some $i\geq 0$.  Due to our assumption on $g$
in this lemma and since $g(r,0)\equiv 0$, we must then have $g(r,u_i) \geq g(r,\widetilde{u}_i) \geq 0$ in $\Theta_T'$.
Next, we apply part (b) of Lemma \ref{pos} to find $\mathscr{L} g(r,u_i) \geq \mathscr{L}g(r,\widetilde{u}_i) \geq 0$ in
$\Theta_T'$.  Proceeding with the iteration argument of the previous lemma, we are thus able to establish
\begin{equation*}
u_{i+1} \equiv u_0 + \mathscr{L}g(r,u_i) \geq \widetilde{u}_0 + \mathscr{L}g(r,\widetilde{u}_i) \equiv
\widetilde{u}_{i+1} \geq 0 \quad \text{in\, $\Theta_T'$.}
\end{equation*}
Since this inequality holds for any $i\geq 0$ by induction, we deduce that $u\geq \widetilde{u} \geq 0$ in $\Theta_T'$.
\end{proof}

\begin{theorem}\label{bu1}
Fix some $R>0$.  Suppose $F\colon \R\to \R$ is a $\mathcal{C}^1$ function with $F(0)=0$ and
\begin{equation}\label{Fb}
F'(u)\geq 0,\quad F(u)\geq Au^p \quad\quad\text{on\, $[0,\infty)$}
\end{equation}
for some $A>0$ and $p>1$.  Suppose also that $V\colon (R,\infty)\to \R$ is $\mathcal{C}^1$ and non-positive.  Let
$\psi\colon (R,\infty) \to \R$ be a continuous function with
\begin{equation}\label{psi1}
\psi(r) \geq \e r^{-k-1} \quad \text{on\, $(R, \infty)$}
\end{equation}
for some $\e>0$ and some $0\leq k< 2/(p-1)$. With $g(r,u)= F(u) - V(r)\cdot u$, consider now the unique solution $u$ to
\eqref{nlb} provided by Lemma \ref{ex}.  If $T$ denotes the lifespan of $u$, then
\begin{equation}\label{ubf}
T\leq C\e^{-(p-1)/(2-k(p-1))}
\end{equation}
for some constant $C$ that is independent of $\e$.
\end{theorem}

\begin{proof}
First of all, note that we may decrease the value of $\e$ without loss of generality.  In particular, we may assume $\e$
is so small that $\e^{-(p-1)/(2-k(p-1))} \geq R$. This trivially gives the desired \eqref{ubf} in the case $T\leq R$, so
we now focus on the case $T\geq R$.

\Step{1} We choose a smooth cut-off function $\zeta \in \mathcal{C}_c(\R)$ and look at the solution $\widetilde{u}$ of
\eqref{nlb} with $\psi$ replaced by $\zeta\psi$.  The idea is to choose $\zeta$ in such a way that
\begin{equation}\label{bl1}
u(r,t) \geq \widetilde{u}(r,t) \geq 0 \quad \text{when\, $(r,t)\in \Theta_T'$}
\end{equation}
as well as
\begin{equation}\label{bl2}
\widetilde{u}(r,t) = 0 \quad\text{when\, $(r,t)\notin \Theta_T'$ \,and\, $0\leq t< T$.}
\end{equation}
The region $\Theta_T'$, which is defined by \eqref{Th2}, corresponds to the shaded region in Figure \ref{cm}.  To ensure
the first inequality \eqref{bl1}, we use our Comparison Lemma \ref{cl}.  In our case,
\begin{equation*}
\d_u g(r,u) = F'(u) - V(r) \geq 0 \quad\quad\text{if\, $(r,u)\in (R,\infty)\times [0,\infty)$}
\end{equation*}
by assumption, so the lemma is applicable.  As long as $0\leq \zeta\leq 1$, we have $\psi\geq \psi\zeta\geq 0$ on
$(R,\infty)$, whence \eqref{bl1} follows.  To ensure the second inequality \eqref{bl2}, we require that
\begin{equation}\label{co}
\zeta(r) = \left\{
\begin{array}{ccc}
0 & \text{if} &  r\leq T_* \\
1 & \text{if} & 2T_* \leq r\leq 3T_* \\
0& \text{if}  & r\geq 4T_*
\end{array}\right\}, \quad\quad T_* = \max((\b_n+2) T, 2T+R).
\end{equation}
Given any $(r,t)\notin \Theta_T'$ with $0\leq t< T$, we then have
\begin{equation*}
r+t\leq \max(\b_n t, t+R) + T \leq T_*.
\end{equation*}
Since $\zeta(r)=0$ for $r\leq T_*$, we must thus have $\widetilde{u}(r,t)=0$ by finite speed of propagation.

\begin{figure}[t]
\begin{pspicture}(0,0)(6.5,4.3) %\showgrid
% Positivity region
\pspolygon[linestyle=none, fillstyle=solid, fillcolor=cyan](1,0)(2,1)(3,1.5)(6.5,1.5)(6.5,0)
\rput(3.2,0.7){\tiny $\Theta_T'$}

% Coordinate axes
\psline{->}(6.5,0)
\rput(6.5,-0.2){\tiny $r$}
\psline{->}(0,4)
\rput(-0.2,4){\tiny $t$}

% Solid lines
\psline(1,0)(4,3)
\rput(0.9,-0.2){\tiny $R$} \rput(4.4,3.2){\tiny $r=t+R$}
\psline(0,0)(5,2.5)
\rput(5.6,2.4){\tiny $r=\b_n t$}
\psline(0,1.5)(6.5,1.5)
\rput(-0.2,1.5){\tiny $T$}

% Dashed line
\rput(5.6,-0.2){\tiny $T_*$}
\psline[linestyle=dashed, dash=4pt 2pt](2.5,3)(5.5,0)
\rput(2.2,3.2){\tiny $r+t= T_*$}
\end{pspicture}
\caption{Our comparison lemma applies only in the shaded region $\Theta_T'$.}\label{cm}
\end{figure}
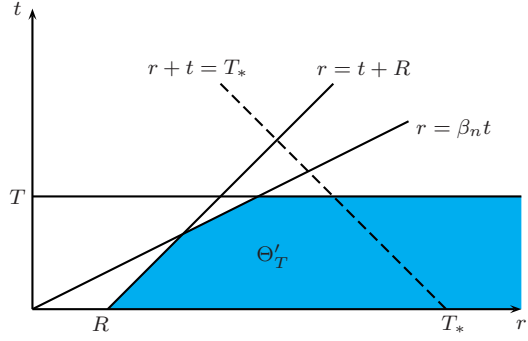

\Step{2} We show that $\widetilde{u}$ becomes infinite before time $T$, unless the estimate \eqref{ubf} holds.  In view
of \eqref{bl1} and \eqref{bl2}, this automatically implies the exact same result for $u$.

As we have already noted, $\widetilde{u}$ vanishes when $r+t\leq T_*$, namely to the left of the dashed line in Figure
\ref{cm}. Although Lemma \ref{ex} does not define $\widetilde{u}$ in the unshaded region to the left of this line, we may
extend $\widetilde{u}$ to be identically zero there.  Thus, we obtain a non-negative solution $\widetilde{u}$ to the
nonlinear problem
\begin{equation}\label{nlb2}
\left\{
\begin{array}{rll}
\d_t^2 \widetilde{u} -\d_r^2 \widetilde{u} - \dfrac{n-1}{r} \cdot \d_r \widetilde{u}= F(\widetilde{u}\,) - V(r)\cdot
\widetilde{u} \quad\quad &\text{in\, $\Omega_T= \R_+\times (0,T)$}\\
\widetilde{u}(r,0)= 0; \quad \d_t \widetilde{u}(r,0)= \psi(r)\zeta(r) \quad\quad &\text{in\, $\R_+$.}
\end{array}\right.
\end{equation}
By \eqref{co} and finite speed of propagation, one has $\widetilde{u}(r,t)=0$, unless $T_*-t \leq r \leq 4T_*+t$. Now,
consider the function
\begin{equation}\label{ft}
f(t) = \int_0^\infty \widetilde{u}(r,t) \cdot r^{n-1} \:dr = \int_{T_\ast -t}^{4T_*+t} \widetilde{u}(r,t) \cdot r^{n-1}
\:dr.
\end{equation}
This function is non-negative on $[0,T)$ and can be easily seen to satisfy
\begin{equation*}
f''(t) = \int_{T_\ast -t}^{4T_*+t} F(\widetilde{u}(r,t)) \cdot r^{n-1} \:dr - \int_{T_\ast -t}^{4T_*+t} V(r)\cdot
\widetilde{u}(r,t) \cdot r^{n-1} \:dr.
\end{equation*}
Since $T_*\geq 2T+R$ by \eqref{co}, one has $r\geq T_*-t \geq R$ within the region of integration, so the potential term
$V(r)$ is non-positive by assumption.  Using our assumption on $F$ together with H\"older's inequality, one then easily
finds
\begin{align}\label{ODE}
f''(t) \geq C(n,p)\cdot (4T_* + T)^{-n(p-1)} \cdot f(t)^p.
\end{align}
This makes $f(t)$ convex, while our initial conditions \eqref{nlb2}  give
\begin{equation*}
f(0)= 0, \quad f'(0) = \int_{T_*}^{4T_*} \zeta(r)\psi(r) \cdot r^{n-1}\:dr.
\end{equation*}
Recalling the definition \eqref{co} of $\zeta$ and our assumption \eqref{psi1} on $\psi$, we thus arrive at
\begin{equation}\label{one}
f(t) \geq t\int_{2T_*}^{3T_*} \psi(r) \cdot r^{n-1}\:dr \geq C(n,k) \cdot \e tT_*^{n-k-1}.
\end{equation}

As we mentioned in the beginning of the proof, we need only treat the case $T\geq R$.  For this case, the constant $T_*=
\max((\b_n+2)T, 2T+R)$ of \eqref{co} is equivalent to $T$ itself, so the analysis simplifies to some extent. Since
$f(0)=0$, an integration of \eqref{ODE} leads to
\begin{equation}\label{ODE2}
f'(t) \geq CT^{-n(p-1)/2} \cdot f(t)^{(p+1)/2}.
\end{equation}
Further integrating on $[T/2,T)$ and using the fact that $p>1$, one now arrives at
\begin{equation*}
f(T/2)^{-(p-1)/2} - CT^{1-n(p-1)/2} \geq \lim_{t\to T} f(t)^{-(p-1)/2}.
\end{equation*}
If the left hand side happens to be negative, then $f(t)\to \infty$ as $t\to T$.  Otherwise, we get
\begin{equation*}
f(T/2) \leq CT^{n- 2/(p-1)}
\end{equation*}
and we also have $f(T/2)\geq C\e T^{n-k}$ by \eqref{one}, so we find that $\e \leq CT^{k-2/(p-1)}$. In view of our
assumption that $k<2/(p-1)$, the desired estimate \eqref{ubf} follows trivially.
\end{proof}

\section{Blow-up due to the Potential}\label{BP}
Our main goal in this section is to give the

\begin{proof_of}{Theorem \ref{bpg}}
We are given a positive and exponentially decaying function $\chi(x)$ which satisfies $(-\Delta + V)\chi = -\la\chi$ for
some $\la>0$, as well as a function $u$ of compact support which is of finite energy and satisfies
\begin{equation*}
\d_t^2 u + (-\Delta + V)u = A|u|^p
\end{equation*}
for some $A>0$ and $p>1$.  Using a standard limiting argument then, we find that
\begin{equation*}
\frac{d^2}{dt^2} \int_{\R^n} \chi u\:dx  - \la \int_{\R^n} \chi u\:dx = A\int_{\R^n} \chi|u|^p \:dx.
\end{equation*}
Since we also have
\begin{equation*}
\int_{\R^n} \chi |u|^p \:dx \geq \left( \int_{\R^n} \chi \:dx \right)^{-(p-1)} \cdot \left| \int_{\R^n} \chi u \:dx
\right|^p
\end{equation*}
by H\"older's inequality, we may now combine the last two equations to arrive at
\begin{equation}\label{ineq}
f''(t) - \la f(t) \geq C(A,\chi,p)\cdot |f(t)|^p, \quad\quad f(t) = \int_{\R^n} \chi u\:dx.
\end{equation}
Our positivity assumptions on the initial data \eqref{pde} are merely imposed to ensure
\begin{equation*}
f(0) = \int_{\R^n} \phi\chi\:dx \geq 0, \quad\quad f'(0) = \int_{\R^n} \psi\chi\:dx > 0.
\end{equation*}
Let $[0,T_0)$ be the maximal interval on which $f'>0$.  Since $f(t)\geq f(0)\geq 0$ on this interval, one also has
$f''(t)\geq 0$ by \eqref{ineq}.  This implies $f'(T_0)\geq f'(0)>0$, whence $f'$ cannot really vanish at $T_0$.  In other
words, $f'$ remains positive as long as it exists, and the same is also true for $f$.  An immediate consequence of
\eqref{ineq} is then
\begin{equation*}
f''(t) \geq Cf(t)^p.
\end{equation*}
Since $p>1$, this ordinary differential inequality implies that $f(t)$ becomes infinite in finite time.  Given any $1\leq
q\leq \infty$, however, we also have
\begin{equation*}
f(t) = \int_{\R^n} \chi(x) u(x,t) \:dx \leq ||u(\cdot\,,t)||_{L^q(\R^n)} \cdot ||\chi||_{L^{q/(q-1)}(\R^n)}
\end{equation*}
by H\"older's inequality, whence $||u(\cdot\,,t)||_{L^q(\R^n)}$ becomes infinite in finite time as well.
\end{proof_of}

\begin{ex}
Let $n\geq 1$ and suppose $V\colon \R^n\to \R$ is continuous with
\begin{equation}\label{bd}
V(x) \leq -V_0 |x|^{-2+\de} \quad\quad \text{if \, $|x|\geq R$}
\end{equation}
for some fixed positive constants $V_0,R$ and $\de$.  Then, $-\Delta + V$ has a negative eigenvalue.  See, for instance,
page 87 in \cite{RS4}.
\end{ex}

\begin{ex}
Let $n=2$ and suppose both $|V(x)|^{1+\de}$ and $|V(x)| \,(1+|x|)^\de$ are integrable for some $\de>0$.   Then, $-\Delta
+ aV$ has a negative eigenvalue for all small positive $a$ if and only if the condition $\int V(x)\:dx < 0$ holds. Thus,
a negative eigenvalue may emerge even for rapidly decaying potentials.  In fact, a similar result is true when $n=1$ as
well; see \cite{Sim2}.
\end{ex}

\begin{ex}
Let $n\geq 3$ and suppose $V\in L^{n/2}(\R^n) + L^\infty (\R^n)$ vanishes at infinity.  If the operator $-\Delta + V$ has
a negative eigenvalue, then there exists a least such eigenvalue which comes with a positive eigenfunction. See chapter
11 in \cite{LL} for more details and the similar hypotheses that one needs when $n=1,2$.
\end{ex}

As for the exponential decay of eigenfunctions, we refer the reader to the survey \cite{Sim1}.

\section{The Riemann Operator in the Radial Case}\label{app}
The proofs of the following three lemmas can be found in \cite{Ka1}.

\begin{lemma}[The Riemann operator]\label{hs}
Letting $z(\la,r,t)$ be the rational function
\begin{equation}\label{z}
z(\la,r,t) = \frac{\la^2+ r^2-t^2}{2r\la} \:,
\end{equation}
we define the Riemann operator $L$ as follows.  When $n$ is odd, we set
\begin{equation}\label{Lo}
[Lf](r,t) = \frac{1}{2r^{(n-1)/2}} \: \int_{|t-r|}^{t+r} \la^{(n-1)/2} f(\la) \cdot P_m(z(\la,r,t))\:d\la,
\end{equation}
where $m= (n-3)/2$ and $P_m$ is the $m$th Legendre polynomial.  When $n$ is even, we set
\begin{align}\label{Le}
[Lf](r,t) &= \frac{1}{2r^{(n-1)/2}} \: \int_{|t-r|}^{t+r} \la^{(n-1)/2} f(\la) \cdot U_m(z(\la,r,t)) \,d\la \notag\\
&\quad + \frac{1}{2r^{(n-1)/2}} \: \int_0^{\max (t-r,0)} \la^{(n-1)/2} f(\la) \cdot U_m(z(\la,r,t)) \,d\la,
\end{align}
where
\begin{equation}\label{U}
U_m(z) = \frac{\sqrt 2}{\pi} \int_{\max(z,-1)}^1 \frac{1}{\sqrt{\sigma - z}} \cdot \frac{T_m(\sigma)}{\sqrt{1- \sigma^2}}
\:\:d\sigma
\end{equation}
with $m= (n-2)/2$ and $T_m$ the $m$th Tchebyshev polynomial.  A solution to the homogeneous wave equation \eqref{he} is
then provided by the formula
\begin{equation}\label{u0}
u_0(r,t) = [L\psi](r,t) + \d_t [L\phi](r,t).
\end{equation}
Moreover, this solution belongs to $\mathcal{C}^1(\R_+^2)$ as long as $\psi\in \mathcal{C}^1(\R_+)$ and $\phi\in
\mathcal{C}^2(\R_+)$.
\end{lemma}

\begin{lemma}\label{Lex}
Let $n\geq 4$ be an integer and define $a,m$ by \eqref{am}. Suppose $f\in \mathcal{C}^1(\R_+)$ satisfies the singularity
condition
\begin{equation*}
f(\la) = O\left( \la^{-2m-2+\de} \right) \quad \text{as $\la \rightarrow 0$}
\end{equation*}
for some fixed $\de>0$.  When $D= (\d_r,\d_t)$ and $t\leq 2r$, the Riemann operator is then such that
\begin{align*}
|D^\b [Lf](r,t)| &\leq C_1(n) \cdot r^{-m-a} \int_{|t-r|}^{t+r} \frac{\la^{m-|\b|}}{(r-t+\la)^{1-a}} \cdot \sum_{s=0}^j
\la^{s+1} \,|f^{(s)}(\la)| \:d\la \notag \\
&\quad + \frac{C_2(n)\cdot r^{-m-a}}{(t-r)^{m+|\b|}} \:\int_0^{\max(t-r,0)} \frac{\la^{2m}}{(t-r-\la)^{1-a}} \cdot
\sum_{s=0}^j \la^{s+1} \,|f^{(s)}(\la)| \:d\la \notag\\
&\quad + C_1(n)\cdot r^{-m-a} \: \sum_{s=0}^{j-1} \left[ \la^{m+a-|\b|+s+1} \,|f^{(s)}(\la)| \right]_{\la= |t\pm r|}
\end{align*}
for each $|\b|\leq j\leq 1$.  Moreover, one has $C_2(n)= 0$ when $n$ is odd.
\end{lemma}

\begin{lemma}\label{Lin}
When $t\geq 2r$, the assumptions of the previous lemma imply
\begin{align*}
|D^\b [Lf](r,t)| &\leq C_1'(n) \cdot r^{j-|\b|-m-a} \int_{t-r}^{t+r} \frac{\la^{m-j}}{(r-t+\la)^{1-a}} \cdot \sum_{s=0}^j
\la^{s+1} \,|f^{(s)}(\la)| \:d\la \notag \\
&\quad + \frac{C_2'(n)\cdot r^{j-|\b|-m}}{(t-r)^{j+m+a}} \: \int_0^{t-r} \frac{\la^{2m}}{(t-r-\la)^{1-a}} \cdot
\sum_{s=0}^j \la^{s+1} \,|f^{(s)} (\la)| \:d\la
\end{align*}
for each $|\b|\leq j\leq 1$.  Moreover, one has $C_2'(n)=0$ when $n$ is odd.
\end{lemma}

\section*{Acknowledgement}
The author would like to express his gratitude to Walter A.~Strauss, under whose direction this research was conducted as
part of the author's doctoral dissertation at Brown University.

\end{document}